\documentclass[final]{siamltex}

\def\D{\Delta}
\def\la{\lambda}

\def\e{{\varepsilon}}

\usepackage{graphicx}
\usepackage{subcaption}
\usepackage{amsmath}
\usepackage{booktabs} 
\usepackage{multirow}
\usepackage{array}
\usepackage{slashed}
\usepackage{amsmath,bm}
\usepackage{appendix}
\usepackage{mathrsfs,psfrag,eepic,epsfig}
\usepackage{fancyhdr}
\usepackage{graphicx}
\usepackage{makecell}
\usepackage{multirow}
\usepackage{float}
\usepackage{amsmath,amssymb}
\usepackage{algorithm}
\usepackage{algorithmic}
\usepackage{tikz}
\usepackage{epstopdf}
\usepackage{appendix}

\usepackage[mathlines]{lineno}

\newcommand{\tabincell}[2]{\begin{tabular}{@{}#1@{}}#2\end{tabular}}

\newtheorem{example}{Example}
\newtheorem{remark}{Remark}[section]

\title{First-order Perturbation Theory of Trust-Region Subproblem}

\author{Bo Feng\thanks{School of Mathematics, China University of Mining and Technology, 221116, Jiangsu, P.R. China. E-mail: {\tt bofeng@cumt.edu.cn}. This  author is supported by the Fundamental Research Funds for the Central Universities under grant 2022XSCX07.}
       \and Gang Wu\thanks{Corresponding author. School of Mathematics,
China University of Mining and Technology, Xuzhou 221116, Jiangsu, P.R. China.
E-mail: {\tt gangwu@cumt.edu.cn}.
}
}

\begin{document}

\maketitle

\begin{abstract}
Trust-region subproblem (TRS) is an important problem arising in many applications such as numerical optimization, Tikhonov regularization of ill-posed problems, and constrained eigenvalue problems. In recent decades, extensive works focus on how to solve the trust-region subproblem efficiently. To the best of our knowledge, there are few results on perturbation analysis of the trust-region subproblem. In order to fill in this gap, we focus on first-order perturbation theory of the trust-region subproblem. The main contributions of this paper are three-fold. First, suppose that the TRS is in \emph{easy case}, we give a sufficient condition under which the perturbed TRS is still in \emph{easy case}. Second, with the help of the structure of the TRS and the classical eigenproblem perturbation theory, we perform first-order perturbation analysis on the Lagrange multiplier and the solution of the TRS, and define their condition numbers. Third, we point out that the solution and the Lagrange multiplier could be well-conditioned even if TRS is in {\it nearly hard case}. The established results are computable, and are helpful to evaluate ill-conditioning of the TRS problem beforehand.
Numerical experiments show the sharpness of the established bounds and the effectiveness of the proposed strategies.
\end{abstract}


\section{Introduction}
 In this paper, we are interested in the first-order perturbation theory of the trust-region subproblem (TRS) \cite{9}, \cite[Chap. 4]{7}:
\begin{equation}\label{1}
\min_{\|\bm x\|_2\leq\D} \left\{{f}(\bm x)=\frac{1}{2}\bm x^T{A}\bm x+\bm x^T{\bm g}\right\},
\end{equation}
where $A\in \mathbb{R}^{n\times n}$ is a symmetric matrix, $\bm 0\neq \bm g\in \mathbb{R}^n$, and $\D>0$.
TRS arises in significant applications such as the regularization or smoothing of discrete forms of ill-posed
problems \cite{M1,M2,M3}, and the trust-region globalization strategy used to force convergence in
optimization methods \cite{M1}. It also stems from graph partitioning problems \cite{W.W.1} and the Levenberg--Marquardt algorithm for solving nonlinear least squares problems \cite{7}. Moreover, solving TRS is a key step in trust-region methods for dealing with general nonlinear optimization problems \cite{9,7}.

A global solution to the TRS \eqref{1} is characterized as follows.
\begin{theorem}\cite{9,28}\label{7.59}
The vector $\bm x_{*}$ is a global optimal solution of the trust-region problem \eqref{1} if and only if $\|\bm x_{*}\|\leq \Delta$ and there exists Lagrange
multiplier $\lambda_{*} \geq 0$ such that
 \begin{equation}\label{eq2100}
(A+\lambda_{*}I){\bm x_{*}}=-\bm g,~\lambda_{*}(\Delta-\|\bm x_{*}\|)=0~~{ and}~~
A+\lambda_{*}I \succcurlyeq \bm O.
\end{equation}
\end{theorem}
Let the eigendecomposition of ${A}$ be
\begin{align*}
{A}=\begin{pmatrix}{U}_1&{U}_2\end{pmatrix}
\begin{pmatrix}
{\Lambda}_1 &  \\
          & {\Lambda}_2
\end{pmatrix}
\begin{pmatrix}
{U}_1^T\\
{U}_2^T
\end{pmatrix}={U}_1{\Lambda}_1{U}_1^T+
{U}_2{\Lambda}_2{U}_2^T,
\end{align*}
where $({U}_1 ~ {U}_2)=[\bm u_1,\bm u_2,\ldots,\bm u_n]\in \mathbb{R}^{n\times n}$ is a unitary matrix, and
$$
{\Lambda}_1=diag({\alpha}_1,\alpha_2\ldots,{\alpha}_{n-s})
,\quad
{\Lambda}_2=diag({\alpha}_{n-s+1},\ldots,{\alpha}_{n}),
$$
with
$
\alpha_1\geq\alpha_2\geq\cdots\geq\alpha_{n-s}>\alpha_{n-s+1}=\cdots=\alpha_n
$
being eigenvalues of $A$. Here $s\geq 1$ is the multiplicity of $\alpha_n$.
Indeed, there are two situations for the TRS \eqref{1}  \cite{W.W,28,7}:

${\bf Case ~1}.$  $\emph{Easy case}$:
$\la_{*}>-\alpha_n$ and $\la_{*}\geq 0$.
In this case, $A+\la_{*}I\succ\bm O$, the solution $\bm x_{*}$ for TRS \eqref{1} is unique and
 $
 \bm x_{*} = -(A+\la_{*}I)^{-1}\bm g.
 $

 ${\bf Case~2 }.$ $\emph{Hard case}$: $\la_{*}=-\alpha_n$.
 {In this case, we have \cite{W.W,4}}
 \begin{equation}\label{eq1.3}
 \cos\angle(\bm g,\mathcal{ U}_2)=0
 ~~{\rm and}~~\|(A-\alpha_nI)^\dag\bm g\|\leq \D,
 \end{equation}
where  $\mathcal{U}_2= \mathcal{R}(U_2)$ is the eigenspace associated with the smallest eigenvalue  $\alpha_n$.

In particular, if $\la_{*}$ is {\it close to} $-\alpha_n$, i.e., $\la_{*}+\alpha_n\rightarrow 0$, we call that TRS is in  {\it nearly hard case} \cite{M1}, which can be viewed as a ``special" \emph{easy case}. In \cite[p.278]{3}, it was shown that mathematically the \emph{hard case} represents only a set of TRS instances
of measure zero, it can happen for matrices with special structures, and numerically
there are \emph{nearly hard case}. Therefore, we pay special attention to the \emph{easy case} in this paper.

In recent decades, extensive methods have been proposed for solving medium sized or large-scale TRS \cite{3,10,28,24,M1,M2,T.S,P.L.,Z.L}, and theoretical results were established for convergence theory or error analysis of these methods \cite{9,BF,G.S.,5,L,28,C.sec,4}.
However, to the best of our knowledge, there are few results on perturbation theory of TRS \eqref{1}. To fill-in this gap,
we try to develop first-order perturbation theory for TRS \eqref{1} in this paper.

Recall that both $\bm g$ and $A$ are real and $A$ is symmetric in the trust-region subproblem.
Let ${A}(\varepsilon)=A+\varepsilon E$ be the perturbed matrix, and ${\bm g}(\varepsilon) = \bm g+\varepsilon \bm e$ be the perturbed vector. Then the {\it perturbed} TRS problem can be described as
\begin{equation}\label{2}
\min_{\|\bm s\|_2\leq\D}\left\{ {f}_{\varepsilon}(\bm s)=\frac{1}{2}\bm s^T{A}(\varepsilon)\bm s+\bm s^T{\bm g(\varepsilon)}\right\},
\end{equation}
where $E\in \mathbb{R}^{n\times n}$, $\bm e\in \mathbb{R}^n$, $\varepsilon\in \mathbb{R}$, and ${A}(\varepsilon)\in \mathbb{R}^{n\times n}$ is symmetric.
Let ${\bm x}_{*}$ be the solution of TRS \eqref{1} and
${\la}_{*}$ be the Lagrange multiplier of TRS \eqref{1}, respectively.
In terms of Theorem \ref{7.59}, there are three cases altogether for the solution $\bm x_*$ and the Lagrange multiplier $\la_*$:
\begin{center}
{\it {\rm (I)} $\la_*>0$,~ {\rm (II)} $\la_*=0$~{\rm with}~ $\|\bm x_*\|<\D$, ~~{\rm and}~~ {\rm (III)} $\la_*=0$~{\rm with}~$\|\bm x_*\|=\D$.}
\end{center}
Thus, we focus on the perturbation theory of TRS \eqref{1} for the above three cases.

The contributions of this work are as follows: First, suppose that the TRS is in \emph{easy case},
we give a sufficient condition under which the perturbed TRS is still in \emph{easy case}. Second, with the help of the structure of the matrix problem arising in TRS \eqref{1}, we perform first-order perturbation analysis on the Lagrange multiplier $\la_*$ and the solution $\bm x_*$ of the TRS \eqref{1}, and define condition numbers for them. Third, we point out that the solution and the Lagrange multiplier could be well-conditioned even if TRS \eqref{1} is in {\it nearly hard case}. The results are illustrated by some examples.

This paper is organized as follows. Assume that the TRS \eqref{1} is in \emph{easy case}, in Section 2, we derive a condition under which the TRS \eqref{1} is still in \emph{easy case} after perturbation. In Section 3, we perform first-order perturbation analysis on the multiplier $\lambda_{*}$ and the TRS solution $\bm x_{*}$, and define condition numbers for them. Examples are given to show the sharpness of the established results. In Section 4, we perform some numerical experiments to show that our results are computable and are practical in use for large-scale TRS. Some concluding remarks are given in Section 5.

In this paper, $A\succcurlyeq\bm O~(A\succ\bm O)$ implies that $A$ is symmetric semi-positive definite (positive definite). Let ${\bm x}_{*}(\varepsilon)$ be the solution, and ${\la}_{*}(\varepsilon)$ be the Lagrange multiplier for the perturbed TRS \eqref{2}, respectively. We denote by $(\cdot)^T$ the transpose of a matrix or vector, by $(\cdot)^\dag$ the  Moore-Penrose inverse of a matrix, by $\|\cdot\|$ the Euclidean norm of a
matrix or vector and by $\mathfrak{Re}(\cdot)$ and $\mathfrak{Im}(\cdot)$ the real and imaginary parts of a complex number, respectively.
Let $\mathcal{W}$ be a linear subspace of $\mathbb{R}^n$, and let $W$ be an orthonormal basis of $\mathcal{W}$.
The cosine of the angle between a nonzero vector $\bm p$ and the subspace $\mathcal{W}$ is defines as \cite{Y.S}
 \begin{equation*}
 \cos\angle(\bm p, \mathcal{W})=\frac{\|W^T\bm p\|}{\|\bm p\|}.
 \end{equation*}
Let $\bm 0$, $\bm O$ and $ I$ be the zero
vector, zero matrix and identity matrix, respectively, whose sizes are clear from the context.

\section{When the TRS \eqref{1} is Still in Easy Case After Perturbation}

Suppose that TRS \eqref{1} is in \emph{easy case}.
In this section, we establish a upper bound on $\e$, such that the perturbed TRS \eqref{2} is still in \emph{easy case}. We first need the following two lemmas.
\begin{lemma}\cite[Corollary 7.4.9.3]{Horn}\label{Lem17.57}
Let $\alpha_1\geq\alpha_2\geq\cdots\geq\alpha_n$ and $\widetilde{\alpha}_1\geq\widetilde{\alpha}_2\geq\cdots\geq\widetilde{\alpha}_n$ be eigenvalues of $A$ and ${A}(\varepsilon)$, respectively. Then
$$
\left\|diag{\big(}\alpha_1,\ldots,\alpha_n{\big)}
-diag{\big(}\widetilde{\alpha}_1,\ldots,\widetilde{\alpha}_n{\big)}\right\|\leq |\varepsilon|\|E\|.
$$
\end{lemma}
\begin{lemma}\cite[Theorem 5.1]{Wdin}\label{le1}
Assume that $rank(B)=rank(\widetilde{B})$ and $\kappa(B)\widetilde{\e}<1$, where $\kappa(B)=\|B\|\|B^\dag\|$ and $\widetilde{\e}=\frac{\|B-\widetilde{B}\|}{\|B\|}$. Then
\begin{align*}
\|B^\dag\bm b-\widetilde{B}^\dag\widetilde{\bm b}\|\leq\frac{\kappa(B)}{\left(1-\widetilde{\e}\right)\|B\|}
&\cdot\left(\widetilde{\e}\|B^\dag\bm b\|\|B\|+\|\bm b-\widetilde{\bm b}\| +\widetilde{\e}\kappa(B)\|(I-BB^\dag)\bm b\|\right)\\
&+\widetilde{\e}\|B\|\|(BB^T)^\dag\bm b\|
\end{align*}
\end{lemma}
We are ready to prove the main theorem in this section.
\begin{theorem}\label{thm2119}
Suppose that TRS \eqref{1} is in easy case.
Let
\begin{align*}
\eta_1=&\frac{1}{2}\min\Bigg\{\frac{(\alpha_{n-s}-\alpha_{n} )\big(\|(A-\alpha_nI)^{\dag}\bm g\|-\D\big)}{\frac{3}{2}\|(A-\alpha_nI)^{\dag}\bm g\|+1+\frac{\cos\angle(\bm g,\mathcal{U}_2)\|\bm g\|}{\alpha_{n-s}-\alpha_n}},~
\alpha_{n-s}-\alpha_n\Bigg\},\\
\eta_2=&\frac{1}{2}\min\left\{\frac{\|U_2U_2^T\bm g\|(\alpha_{n-s}-\alpha_n)}{4\|\bm g\|+(\alpha_{n-s}-\alpha_n)},\frac{\|\bm g\|\cos\angle(\bm g, \mathcal{U}_2)}{2\D},~\alpha_{n-s}-\alpha_n\right\}.
\end{align*}
Then the perturbed TRS \eqref{2} is  also  in  easy  case if
\begin{equation}\label{eq}
0<\max\{\|\e E\|, \|\e\bm e\|\}<\max\{ \eta_1,\eta_2\}.
\end{equation}
Specifically, if $A \succ\bm O$, then TRS \eqref{2} is also in \emph{easy case} if $0<\|\e E\|<\alpha_n$.
\end{theorem}
\begin{proof}
Let the eigendecomposition of ${A}(\varepsilon)$ be
\begin{align}\label{eq1813}
{A}(\varepsilon)=\begin{pmatrix}\widetilde{U}_1 & \widetilde{U}_2\end{pmatrix}
\begin{pmatrix}
\widetilde{\Lambda}_1 &  \\
          & \widetilde{\Lambda}_2
\end{pmatrix}
\begin{pmatrix}
\widetilde{U}_1^T\\
\widetilde{U}_2^T
\end{pmatrix}=\widetilde{U}_1\widetilde{\Lambda}_1\widetilde{U}_1^T+
\widetilde{U}_2\widetilde{\Lambda}_2\widetilde{U}_2^T,
\end{align}
where $\big(\widetilde{U}_1 ~~ \widetilde{U}_2\big)\in \mathbb{R}^{n\times n}$ is a unitary matrix, $\widetilde{\Lambda}_1=diag(\widetilde{\alpha}_1,\ldots,
\widetilde{\alpha}_{n-s})$ and $\widetilde{\Lambda}_2=diag(\widetilde{\alpha}_{n-s+1},\ldots,
\widetilde{\alpha}_{n})$, with $\widetilde{\alpha}_1\geq\widetilde{\alpha}_2\geq\cdots\geq\widetilde{\alpha}_n$.
It follows that
\begin{small}
\begin{align}
{\big\|}({A}(\varepsilon)-\widetilde{\alpha}_nI)^{\dag}{\bm g}(\varepsilon){\big\|}^2
=&
{\big\|}\widetilde{U}_1(\widetilde{\Lambda}_1
-\widetilde{\alpha}_nI)^\dag\widetilde{U}_1^T{\bm g}(\varepsilon)+
\widetilde{U}_2(\widetilde{\Lambda}_2-\widetilde{\alpha}_n I)^{\dag}\widetilde{U}_2^T{\bm g}(\varepsilon){\big\|}^2\nonumber\\
=&{\big\|}\big[\widetilde{U}_1(\widetilde{\Lambda}_1
-\widetilde{\alpha}_nI)\widetilde{U}_1^T\big]^\dag{\bm g}(\varepsilon)\big\|^2+\big\|
\big[\widetilde{U}_2(\widetilde{\Lambda}_2-\widetilde{\alpha}_n I)\widetilde{U}_2^T\big]^{\dag}{\bm g}(\varepsilon){\big\|}^2.\label{eq2214}
\end{align}
\end{small}
Recall from \eqref{eq1.3} that TRS \eqref{1} is in {\it easy case} if
$$
\|(A-\alpha_n I)^\dag\bm g\|>\D~~{\rm or}~~\cos\angle(\bm g, \mathcal{U}_2)>0.
$$
Thus, it is only necessary to consider the following two cases: $\|(A-\alpha_n I)^\dag\bm g\|>\D$ and $\cos\angle(\bm g, \mathcal{U}_2)>0$.
Let us discuss them in more details.

(i) First, we consider the case of $\|(A-\alpha_n I)^\dag\bm g\|>\D$.
Suppose that
$$0<|\e|\cdot\max\{\|E\|,\|\bm e\|\}<\frac{\alpha_{n-s}-\alpha_n}{2}.$$
By Lemma \ref{Lem17.57},
\begin{align}
\widetilde{\alpha}_{n-s}\!-\!\widetilde{\alpha}_{n}
=\widetilde{\alpha}_{n-s}\!-\!{\alpha}_{n-s}
+{\alpha}_{n-s}\!-\!{\alpha}_{n}+{\alpha}_{n}\!-\!\widetilde{\alpha}_{n}
\geq({\alpha}_{n-s}\!-\!{\alpha}_{n})\!-\!2|\varepsilon|\|E\|
>0.\label{eq1535}
\end{align}
Thus,
$\widetilde{\Lambda}_1-\widetilde{\alpha}_n I\succ\bm O$ and $rank(\widetilde{U}_1(\widetilde{\Lambda}_1-\widetilde{\alpha}_n I)\widetilde{U}^T_1)=n-s$.
Note that $I-(A-\alpha_nI)(A-\alpha_nI)^{\dag}=I-U_1U_1^T=U_2U_2^T$. It follows from Lemma \ref{le1} that
{\small
\begin{align}
     &{\Big |}{\big\|}(A-\alpha_nI)^{\dag}\bm g{\big\|}-{\big\|}[\widetilde{U}_1(\widetilde{\Lambda}_1
     -\widetilde{\alpha}_n I)\widetilde{U}^T_1]^\dagger{\bm g}(\varepsilon){\big\|} {\Big |}\nonumber\\
\leq&\frac{2|\e|}{\alpha_{n-s}-\alpha_n}\left(\|(A-\alpha_nI)^{\dag}\bm g\|\|E\|+\|\bm e\|+\frac{\|E\|\|U_2U_2^T\bm g\|}{\alpha_{n-s}-\alpha_n}\right)+|\e|\|E\|
\left\|[(A-\alpha_nI)^2]^{\dag}\bm g\right\|\nonumber\\
\leq&\frac{2|\e|}{\alpha_{n-s}-\alpha_n}\left(\frac{3}{2}\|(A-\alpha_nI)^{\dag}\bm g\|\|E\|+\|\bm e\|+\frac{\|E\|\cos\angle(\bm g,\mathcal{U}_2)\|\bm g\|}{\alpha_{n-s}-\alpha_n}\right)\nonumber\\
\leq&\frac{2|\e|\max\{\|E\|,\|\bm e\|\}}{\alpha_{n-s}-\alpha_n}\left(\frac{3}{2}\|(A-\alpha_nI)^{\dag}\bm g\|+1+\frac{\cos\angle(\bm g,\mathcal{U}_2)\|\bm g\|}{\alpha_{n-s}-\alpha_n}\right).\label{eq11..50}
\end{align}}
Thus, if
\begin{align*}
\max\{\|\e E\|,\|\e\bm e\|\}<\frac{(\alpha_{n-s}-\alpha_{n} )\big(\|(A-\alpha_nI)^{\dag}\bm g\|-\D\big)}{3\|(A-\alpha_nI)^{\dag}\bm g\|+2+\frac{2\cos\angle(\bm g,\mathcal{U}_2)\|\bm g\|}{\alpha_{n-s}-\alpha_n}},
\end{align*}
then
$${\Big |}{\big\|}(A-\alpha_nI)^{\dag}\bm g{\big\|}-{\big\|}[\widetilde{U}_1(\widetilde{\Lambda}_1-\widetilde{\alpha}_n I)\widetilde{U}^T_1]^\dagger{\bm g}(\varepsilon){\big\|} {\Big |} <\|(A-\alpha_nI)^{\dag}\bm g\|-\D,
$$
and thus
${\big\|}[\widetilde{U}_1(\widetilde{\Lambda}_1-\widetilde{\alpha}_n I)\widetilde{U}^T_1]^\dagger{\bm g}(\varepsilon){\big\|}>\D$. By \eqref{eq2214}, ${\big\|}({A}(\varepsilon)-\widetilde{\alpha}_nI)^{\dag}{\bm g}(\varepsilon){\big\|}>\D$, i.e., TRS \eqref{2} is also in \emph{easy case}
when
\begin{small}
\begin{align*}
0<\max\{\|\e E\|,\|\e \bm e\|\}<\eta_1=\min\Bigg\{\frac{(\alpha_{n-s}-\alpha_{n} )\big(\|(A-\alpha_nI)^{\dag}\bm g\|-\D\big)}{3\|(A-\alpha_nI)^{\dag}\bm g\|+2+\frac{2\cos\angle(\bm g,\mathcal{U}_2)\|\bm g\|}{\alpha_{n-s}-\alpha_n}},~
\frac{\alpha_{n-s}-\alpha_n}{2}\Bigg\}.
\end{align*}
\end{small}

(ii) Second, we consider the case of $\cos\angle(\bm g, \mathcal{U}_2)>0$.
It holds that
\begin{align}
\|\widetilde{ U}_2\widetilde{ U}^T_2{\bm g}(\varepsilon) - U_2 U^T_2\bm g\|&\leq
\|\widetilde{U}_2\widetilde{U}^T_2{\bm g}(\varepsilon)-\widetilde{U}_2\widetilde{ U}^T_2\bm g\| +\|\widetilde{U}_2\widetilde{U}^T_2\bm g-U_2U^T_2\bm g\|\nonumber\\
&\leq |\varepsilon|\|\bm e\|+ \|U_2U^T_2-\widetilde{U}_2\widetilde{U}^T_2\|\|\bm g\|.\label{eq1509}
\end{align}
Next, we consider $\|U_2U^T_2-\widetilde{U}_2\widetilde{U}^T_2\|$. Suppose that
\begin{small}
\begin{align}\label{eq22.10}
0<\max\{\|\e E\|,\|\e\bm e\| \}<\frac{\alpha_{n-s}-\alpha_n}{2}.
\end{align}
\end{small}
It follows from Lemma \ref{Lem17.57} that
{\small
\begin{align}
\alpha_{n-s}\!-\!\widetilde{\alpha}_{n-s+1}\!=\!
\alpha_{n-s}\!-\!\alpha_{n}+\alpha_{n-s+1}\!-\!\widetilde{\alpha}_{n-s+1}
\geq\alpha_{n-s}\!-\!\alpha_{n}\!-\!|\varepsilon|\|E\|\!\geq\! \frac{\alpha_{n-s}\!-\!\alpha_{n}}{2}\!>\!0.\label{eq1632}
\end{align}}
If we denote by $\zeta=\frac{1}{2}(\alpha_{n-s}-\widetilde{\alpha}_{n-s+1})$, then
$
\alpha_{n-s}-\zeta>\widetilde{\alpha}_{n-s+1},
$
and
\begin{equation*}
\{\widetilde{\alpha}_{n-s+1},\ldots,\widetilde{\alpha}_{n}\}\subseteq
[\widetilde{\alpha}_{n}, \widetilde{\alpha}_{n-s+1}]\subseteq \mathbb{R}\setminus [\alpha_{n-s}-\zeta,\alpha_{1}+\zeta].
\end{equation*}
Notice that $\{\alpha_{1},\alpha_2,\ldots,\alpha_{n-s}\}\subseteq [\alpha_{n-s},\alpha_{1}]$.
From \cite[p. 251 Theorem 3.6]{G.W.},
\begin{align*}
\|U_2U^T_2-\widetilde{U}_2\widetilde{U}^T_2\|
\leq \frac{\|(A+\varepsilon E)U_1-U_1\Lambda_1\|}{\zeta}
= \frac{|\varepsilon|\|EU_1\|}{\zeta}\overset{\eqref{eq1632}}{\leq} \frac{4|\varepsilon|\|E\|}{\alpha_{n-s}-\alpha_{n}}.
\end{align*}
We see from \eqref{eq1509} that
{\small
\begin{equation}\label{eq21.52}
{\Big|}\|\widetilde{U}_2\widetilde{ U}^T_2\!{\bm g}(\varepsilon)\|-\| U_2 U^T_2\!\bm g\|{\Big |}\!\leq\! |\varepsilon|\|\bm e\|+\frac{4|\varepsilon|\|E\|\|\bm g\|}{\alpha_{n-s}\!-\!\alpha_{n}}\!\leq\!
|\e|\max\{\|E\|,\|\bm e\|\}\left(1\!+\!\frac{4\|\bm g\|}{\alpha_{n-s}\!-\!\alpha_{n}}\right).
\end{equation}}

If $\max\{\|\e E\|,\|\e\bm e\|\}<\frac{\|U_2U_2^T\bm g\|(\alpha_{n-s}-\alpha_n)}{2[4\|\bm g\|+(\alpha_{n-s}-\alpha_n)]}$, we have
$$
|\varepsilon|\|\bm e\|+\frac{4|\varepsilon|\|E\|\|\bm g\|}{\alpha_{n-s}-\alpha_{n}}
<\frac{1}{2}\|U_2U_2^T\bm g\|.
$$
Thus, ${\big \|}\widetilde{ U}_2\widetilde{ U}^T_2{\bm g}(\varepsilon){\big \|}>\frac{1}{2}\|U_2U_2^T\bm g\|>0$. Notice that $0<2\|E\||\e|<\alpha_{n-s}-\alpha_n$, from
 \eqref{eq1535}, we see that $\widetilde{\alpha}_{n-s}>\widetilde{\alpha}_n$. That is, as an eigenvalue of $A(\e)$, the multiplicity $\widetilde{s}$ of $\widetilde{\alpha}_n$ will not be larger than $s$. Let $\widetilde{\mathcal{N}}$ be the eigenspace of $A(\e)$ associated with the smallest eigenvalue $\widetilde{\alpha}_n$.
  On one hand, if $s=\widetilde{s}$, then $\widetilde{\alpha}_n=\widetilde{\alpha}_{n-s+1}$, and $\mathcal{R}(\widetilde{U}_2)=\widetilde{\mathcal{N}}$. Thus, we have from $$\cos\angle(\bm g(\e),\widetilde{\mathcal{N}})=\frac{\|\widetilde{ U}_2\widetilde{ U}^T_2{\bm g}(\varepsilon) \|}{\|\bm g(\e)\|}>0$$ that TRS \eqref{2} is in \emph{easy case}.

 On the other hand, if $\widetilde{s}<s$, then
  $$
 \widetilde{\alpha}_{n-s+1}\geq\cdots
 \geq\widetilde{\alpha}_{n-\widetilde{s}}>
 \widetilde{\alpha}_{n-\widetilde{s}+1}=\cdots=\widetilde{\alpha}_n.
 $$
 Partition $\widetilde{U}_2\in \mathbb{R}^{n\times s}$ as $\widetilde{U}_{2}=\big(\widetilde{U}^{(1)}_2~\widetilde{U}^{(2)}_2\big)$, with $\widetilde{U}^{(1)}_2\in \mathbb{R}^{n\times (s-\widetilde{s})}$ and $\widetilde{U}^{(2)}_2\in \mathbb{R}^{n\times \widetilde{s}}$.
 Here $\widetilde{\mathcal{N}}=\mathcal{R}(\widetilde{U}^{(2)}_2)$, and
 $${\big \|}\widetilde{ U}_2\widetilde{ U}^T_2{\bm g}(\varepsilon){\big \|}^2={\big \|}\widetilde{ U}^T_2{\bm g}(\varepsilon){\big \|}^2=\big \|(\widetilde{U}^{(1)}_2)^T\bm g(\e)\big \|^2+\big \|(\widetilde{U}^{(2)}_2)^T\bm g(\e)\big \|^2.$$
 If $\cos\angle(\bm g(\e),\widetilde{\mathcal{N}})=\frac{\|(\widetilde{U}^{(2)}_2)^T\!\bm g(\e) \|}{\|\bm g(\e)\|}>0$, then TRS \eqref{2} is in \emph{easy case}. If $\cos\angle(\bm g(\e),\widetilde{\mathcal{N}})=\frac{\|(\widetilde{U}^{(2)}_2)^T\!\bm g(\e) \|}{\|\bm g(\e)\|}=0$, we have $\!\big \|(\widetilde{U}^{(1)}_2)^T\!\bm g(\e)\big \|=\!{\big \|}\widetilde{ U}_2\widetilde{ U}_2^T\!{\bm g}(\varepsilon){\big \|}\!>\!\frac{1}{2}\|U_2U^T_2{\bm g}\|\!=\!\frac{\|\bm g\|}{2}\cos\angle(\bm g,\mathcal{U}_2)$. Moreover,
 {\small
 \begin{align*}
 {\Big\|}\big[\widetilde{U}_2(\widetilde{\Lambda}_2\!-\!\widetilde{\alpha}_n I)\widetilde{U}_2^T\big]^{\dag}\!{\bm g}(\varepsilon){\Big\|}\!=\!
 {\Big\|}(\widetilde{\Lambda}^{(1)}_2\!-\!\widetilde{\alpha}_n I)^{-1}(\widetilde{U}^{(1)}_2)^T\!{\bm g}(\varepsilon){\Big\|}
 \!\geq\!
 \frac{{\big\|}(\widetilde{U}^{(1)}_2)^T\!{\bm g}(\varepsilon){\big\|}}{\widetilde{\alpha}_{n-s+1}\!-\!\widetilde{\alpha}_n}
 \!>\!\frac{\|\bm g\|\!\cos\angle(\bm g,\mathcal{U}_2)}{4|\e|\|E\|},
 \end{align*}}
 where $\widetilde{\Lambda}^{(1)}_2=diag(\widetilde{\alpha}_{n-s+1},
 \ldots,\widetilde{\alpha}_{n-\widetilde{s}})$, and the last inequality is from the fact that
 \begin{align*}
 \widetilde{\alpha}_{n-s+1}-\widetilde{\alpha}_n&=
\widetilde{\alpha}_{n-s+1}-{\alpha}_{n-s+1}+
 {\alpha}_{n}-\widetilde{\alpha}_n\overset{{\rm Lem.}~ \ref{Lem17.57}}{\leq} 2|\e|\|E\|.
 \end{align*}
If
$
|\e|\|E\|<\frac{\|\bm g\|\cos\angle(\bm g, \mathcal{U}_2)}{4\D}
$, then
 \begin{align*}
{\big\|}({A}(\varepsilon)-\widetilde{\alpha}_nI)^{\dag}{\bm g}(\varepsilon){\big\|}
\overset{\eqref{eq2214}}{\geq}{\big\|}\widetilde{U}_2(\widetilde{\Lambda}_2-\widetilde{\alpha}_n I)^{\dag}\widetilde{U}_2^T{\bm g}(\varepsilon){\big\|}>\frac{4\|U_2^T\bm g\|}{|\e|\|E\|}> \D.
\end{align*}
Thus,
$
\cos\angle(\bm g(\e), \widetilde{\mathcal{N}})>0~~{\rm or}~~\|(A({\e})-\widetilde{\alpha}_n I)^\dag\bm g(\e)\|-\D>0,
$
provided {\small
$$
0<\max\{\|\e E\|,\|\e \bm e\|\}<\eta_2=\min\left\{\frac{\alpha_{n-s}-\alpha_n}{2},\frac{\|U_2U_2^T\bm g\|(\alpha_{n-s}-\alpha_n)}{8\|\bm g\|+2(\alpha_{n-s}-\alpha_n)},\frac{\|\bm g\|\cos\angle(\bm g, \mathcal{U}_2)}{4\D}\right\}.
$$}
That is,
TRS \eqref{2} is in \emph{easy case} if $\cos\angle(\bm g,\mathcal{U}_2)>0$ and $0<\max\{\|\e E\|,\|\e \bm e\|\}<\eta_2$. In conclusion, TRS \eqref{2} is in \emph{easy case} if $0<\max\{\|\e E\|,\|\e \bm e\|\}<\max\{\eta_1,\eta_2\}$.

(iii) Specially, if $A\succ\bm O$ and $|\varepsilon|\|E\|\!<\!\alpha_n$, then we have from Lemma \ref{Lem17.57} that
$\widetilde{\alpha}_n\!\geq \!\alpha_n\!-\!|\varepsilon|\|E\|>0$.
As a result, ${A}(\varepsilon)\!\succ\!\bm O$, and it follows that TRS \eqref{2} is in \emph{easy case}.
\end{proof}
\section{First-Order Perturbation Theory on the Multiplier $\bm \lambda_{*}$ and the TRS Solution $\bm x_{*}$}
In this section, we focus on the first-order perturbation theory on the multiplier $\lambda_{*}$ and the TRS solution $\bm x_{*}$.
Without loss of generality, we assume that \eqref{eq} is satisfied from now on, such that the perturbed  TRS \eqref{2} is still in \emph{easy case}.
Recall that $\la_*\geq 0$, $\|\bm x_*\|\leq \D$ and $\la_*(\D-\|\bm x_*\|)=0$.
Thus, there are the following three situations altogether:
\begin{center}
{\it {\rm (I)} $\la_*>0$,~ {\rm (II)} $\la_*=0$~{\rm with}~$\|\bm x_*\|=\D$, ~{\rm and}~{\rm (III)} $\la_*=0$~{\rm with}~ $\|\bm x_*\|<\D$.}
\end{center}

It is seen from Theorem \ref{7.59} that $A\succ\bm O$ and
$\bm x_*=A^{-1}\bm g$ in Case (III).
Thus, we have from the classical perturbation theorem for linear systems \cite[section 1.13.2]{2} that
$$
\bm x_*(0)=\lim_{\e\rightarrow 0}\bm x_*(\e)=\bm x_*~~{\rm and}~~
\bm x'_*(0)=\lim_{\e\rightarrow 0}\frac{\bm x_*(\e)-\bm x_*(0)}{\e}=-A^{-1}(E\bm x_*+\bm e).
$$
Moreover, if $\e$ small enough, we have $\|\bm x_*(\e)\|<\D$ and
$\la_*(\e)=\la_*=0$. So we only need to consider Case (I) and (II) in this section.

\subsection{The case of $\bm{\la_{*}>0}$}
In this case, the solution $\bm x_*$ to the TRS \eqref{1} reaches the trust region boundary, i.e.,  $\|\bm x_{*}\|=\D$; see \eqref{eq2100}.
First, we consider the continuity and differentiability of $\la_*(\e)$ at $\e=0$.
It was shown that \eqref{1} can be rewritten as an eigenvalue problem corresponding to
the following $2n$-by-$2n$ matrix \cite{3,GGM,5}:
\begin{equation}\label{eq16.12}
M = \begin{pmatrix}
-A & \frac{\bm g\bm g^T}{\D^2}\\
I  &  -A
\end{pmatrix}\in \mathbb{R}^{2n\times 2n}.
\end{equation}
The following theorem establishes an important
relationship between $(\la_{*},\bm x_{*})$ and the rightmost eigenpair of $M$.
\begin{theorem}\label{16.0200}\cite{3,5}
Let $(\lambda_{*}, \bm x_{*})$ satisfy Theorem \ref{7.59} with $\|\bm x_{*}\| = \D$.
Then the rightmost eigenvalue $\la_R$ of $M$ is real and simple, and $\la_{*} = \la_R$.
  Let $\bm y = (\bm y^T_1,\bm y_2^T )^T$
be the corresponding unit length eigenvector of $M$ with $\bm y_1, \bm y_2 \in \mathbb{R}^{n}$, and suppose that $\bm g^T \bm y_2 \not = 0$. Then the unique TRS solution is
\begin{equation}\label{eq21.24}
\bm x_{*} = -\frac{\D^2}{\bm g^T\bm y_2}\bm y_1=-sign(\bm g^T\bm y_2)\D\cdot\frac{\bm y_1}{\|\bm y_1\|}.
\end{equation}
\end{theorem}
\begin{remark}
It was shown that $\bm g^T\bm y_2=0$ only if TRS is in hard case \cite[Proposition 4.1]{3}.
Moreover, if $\|\bm x_*\|=\D$, then $\la_*$ is simple if and only if TRS \eqref{1} is in easy case \cite[Theorem 4.2]{BF}.
\end{remark}

Since $\la_*$ is the rightmost eigenvalue of $M$ when $\|\bm x_*\|=\D$, based on the first-order eigenvalue perturbation theory \cite{Li,Horn,G.W.},
it seems that one can derive some first-order perturbation results for $\la_*$. However, the key is how to establish refined bounds by exploiting the {\it structure} of the matrix $M$ sufficiently.

To derive the main results, we first need two lemmas.
Theorem \ref{16.0200} indicates that $\la_*\geq\mathfrak{Re}(\la)$ for any $\la\in\la(M)\backslash \{\la_*\}$. The first lemma shows that the inequality holds strictly.
\begin{lemma}\label{lem20.59}
If $\la_*>0$
and TRS \eqref{1} is in easy case, then $\la_*$ is a simple eigenvalue of $M$ and
$\la_*>\mathfrak{Re}(\la)$ for any $\la\in\la(M)\backslash \{\la_*\}$.
\end{lemma}
\begin{proof}
It follows from \cite[Theorem 4.2]{BF} that $\la_*$ is simple.
Recall that $\la_*$ is the rightmost eigenvalue of $M$, and thus $\la_*\geq\mathfrak{Re}(\la)$ holds for any $\la\in\la(M)\backslash\{\la_*\}$.
Assume that $\la_*+c{\bf i}\in \la(M)$, where ${\bf i}^2=-1$ and $c\neq 0$. Notice that $\det(A+(\la_*+c{\bf i})I)\neq 0$. It follows that
\begin{align*}
\det(M-(\la_*+c{\bf i})I)
 &=\det(A+(\la_*+c{\bf i}) I)^2\cdot
 \left(1\!-\!\frac{\bm g^T(A+(\la_*+c{\bf i}) I)^{-2}\bm g}{\D^2} \right)\\
 &=\frac{\det(A+(\la_*+c{\bf i}) I)^2}{\D^2}\cdot
 \left(\D^2- \sum_{i=1}^n\frac{(\bm u_i^T\bm g)^2}{(\alpha_i+\la_*+c{\bf i})^2}\right).
\end{align*}
As $\la_*+\alpha_n>0$ if TRS \eqref{1} is in \emph{easy case} and $c\neq 0$, then $\mathfrak{Im}\left(\sum_{i=1}^n\frac{(\bm u_i^T\bm g)^2}{(\alpha_i+\la_*+c{\bf i})^2}\right)\neq 0$. Consequently,
$\det(M-(\la_*+c{\bf i})I)\neq 0$, which is a contradiction. As a result, $\la_*>\mathfrak{Re}(\la)$ for any $\la\in\la(M)\backslash \{\la_*\}$.
\end{proof}

Motivated by \cite[Theorem 6.3.12]{Horn}, we have the second lemma.
\begin{lemma}\label{19.0333}
If $\la_*>0$. Let
$\bm y=(\bm y_1^T~\bm y_2^T)^T\in \mathbb{R}^{2n}$  be the {\rm(}unit{\rm)} right eigenvector of $M$ corresponding to $\la_{*}$, with $\bm y_1, \bm y_2\in \mathbb{R}^n$. Denote by
 \begin{equation}\label{eq3.6}
 {M(\e)}=\begin{pmatrix}
-{A(\e)}  &  \frac{{\bm g(\e)} {\bm g(\e)}^T}{\D^2}\\
I     & -{A(\e)}
\end{pmatrix}.
 \end{equation}
Then there exists a scalar $\breve{\varrho}>0$, if $0<|\varepsilon|<{\breve{\varrho}}$, there is a simple eigenvalue $\breve{\la}(\varepsilon)$ of $M(\varepsilon)$, such that
\begin{itemize}
\item[\rm(i)]
$
\underset{\e\rightarrow 0}{\lim}\breve{\la}(\e)=\la_*,
$
\item[\rm(ii)]
 $\breve{\la}(\varepsilon)$ is differentiable at $\varepsilon = 0$, and
\begin{eqnarray*}
\breve{\la}'(0)=\frac{{\rm d}\breve{\la}(\varepsilon)}{{\rm d}\varepsilon}{\bigg|_{\varepsilon=0}}=\frac{sign(\bm g^T\bm y_2)\|\bm y_1\|}{\bm y_1^T\bm y_2\D}\cdot\bm y_2^T\left(E\bm x_*+\bm e\right).
\end{eqnarray*}
\end{itemize}
\end{lemma}
\begin{proof}
Notice that
\begin{align}\label{eq1023}
\bm w=(\bm y_2^T
~\bm  y^T_1 )^T
\end{align}
is the (unit) left eigenvector of $M$ associated with $\lambda_*$ \cite[p. 895]{5}, and
$$
M(\e)-M=\e\cdot\begin{pmatrix}
-E &  \frac{\bm e\bm g^T+\bm g\bm e^T}{\D^2}+\e\cdot\frac{\bm e\bm e^T}{\D^2}\\
\bm O  &  -E
\end{pmatrix}\equiv\e F_{\e}.
$$
As $F_{\e}$ is uniformly bounded for $|\e|< 1$ and  $\la_*$ is simple,
we have from \cite[Theorem 6.3.12 (a)]{Horn} that there exists a scalar $\breve{\varrho}>0$, such that if $0<|\varepsilon|<{\breve{\varrho}}$, then there is a simple eigenvalue $\breve{\la}(\varepsilon)$ of $M(\varepsilon)$ such that
$\underset{\e\rightarrow 0}{\lim}\breve{\la}(\e)=\la_*$. Moreover,
 $$
 \left| \frac{\breve{\la}(\e)-\la_*}{\e}- \frac{\bm w^TF_{\e}\bm y}{\bm w^T\bm y} \right|<\mathcal{O}(\e),\quad\quad 0<|\e|<\breve{\varrho}.
 $$
Let $\e\rightarrow 0$, then
\begin{align*}
\lim_{\e\rightarrow 0}\frac{\breve{\la}(\e)-\la_*}{\e}
=&(\bm w^T\bm y)^{-1}\cdot\bm w^T\begin{pmatrix}
-E &  \frac{\bm e\bm g^T+\bm g\bm e^T}{\D^2}\\
\bm O  &  -E
\end{pmatrix}\bm y
=\frac{(\bm g^T\bm y_2)\cdot(\bm y_2^T \bm e)}{\bm y_1^T\bm y_2\D^2}-\frac{\bm y^T_2E\bm y_1 }{\bm y_1^T\bm y_2}\\
=&\frac{(\bm g^T\bm y_2)\cdot(\bm y_2^T \bm e)}{\bm y_1^T\bm y_2\D^2}-\frac{\bm y^T_2E\bm y_1 }{\bm y_1^T\bm y_2}
\overset{\eqref{eq21.24}}{=}\frac{sign(\bm g^T\bm y_2)\|\bm y_1\|\cdot\bm y_2^T\bm e}{(\bm y_1^T\bm y_2)\D}-\frac{\bm y^T_2E\bm y_1}{\bm y_1^T\bm y_2}\\
=&\frac{sign(\bm g^T\bm y_2)\|\bm y_1\|}{\bm y_1^T\bm y_2\D}\cdot\bm y_2^T\left(\bm e-\frac{\D sign(\bm g^T\bm y_2)E\bm y_1}{\|\bm y_1\|}\right)\\
=&\frac{sign(\bm g^T\bm y_2)\|\bm y_1\|}{\bm y_1^T\bm y_2\D}\cdot\bm y_2^T\left(E\bm x_*+\bm e\right),
\end{align*}
where we used $E=E^T$.
\end{proof}
\begin{remark}
We have from Theorem \ref{16.0200} that $\la_*(\e)$ is {\it the rightmost} eigenvalue $\la_R(\e)$ of $M(\e)$ for $\|\bm x_*(\e)\|=\D$. Moreover, the matrix $M(\varepsilon)$ defined in \eqref{eq3.6} is a structured perturbation to $M$.
\end{remark}

Lemma \ref{19.0333} indicates that $\breve{\la}(\e)$ is an eigenvalue of $M(\e)$. For sufficiently small values of $\e$,
we have from the continuity argument that
$\breve{\la}(\e)$ is still the rightmost eigenvalue of $M(\e)$, and all the other eigenvalues have
imaginary parts that are smaller than $\breve{\la}(\e)$. To derive the differentiability of the Lagrange multiplier $\la_*(\e)$ at 0, however,
we have to consider the size of $\e$ such that $\breve{\la}(\e)$ is still the
rightmost eigenvalue of $M(\e)$.

Denote by
$$
\nu=\min_{\la\in \la(M)\backslash \{\la_{*}\}}|\la_{*}-\la|,~~{\rm and}~~\widetilde{\nu}=\min_{\la\in \la(M)\backslash
\{\la_{*}\}}(\la_{*}-\mathfrak{Re}(\la)).
$$
It is seen from Lemma \ref{lem20.59} that $\nu>0$ and $\widetilde{\nu}>0$. We have from \cite[Theorem 6.3.12 (6.3.13a)]{Horn} that there exists a matrix $N$, such that
\begin{equation}\label{eq15.10}
N^{-1}(M+\e F_{\e})N=\begin{pmatrix}
\la_*+\e \frac{\bm w^TF_{\e}\bm y}{\bm w^T\bm y} & \e\bm \zeta^T\\
           \e\bm \beta                &   T+\e W
\end{pmatrix}\in \mathbb{R}^{2n\times 2n},
\end{equation}
where
 $$
 \max\left\{\breve{\varrho}\left|\frac{\bm w^TF_{\e}\bm y}{\bm w^T\bm y}\right|, \breve{\varrho}\|\bm \beta\|_{\infty}, 2\breve{\varrho}\|W\|_{\infty},\|\bm \zeta\|_1\right\}<\frac{\nu}{7},
 $$
and $T\in \mathbb{R}^{(2n-1)\times(2n-1)}$ is an upper triangular matrix with elements $\{t_{ij}\}$'s, and the eigenvalues of $M$ (excluding $\la_*$) on its main diagonal.
Moreover, from \cite[Theorem 6.3.12]{Horn}, we have that the first Ger$\rm \check{s}$gorin disc $G_1$ with respect to the matrix defined in \eqref{eq15.10} is disjoint from the other Ger$\rm \check{s}$gorin discs associated with the rows $2,3,\ldots,2n$, and if $0<|\e|<\breve{\varrho}$, then
\begin{small}
\begin{align}\label{eq1451}
\breve{\la}(\e)\!\in\! G_1\!=\!\left\{\!z:~{\Bigg|}z\!-\!{\bigg(}\la_*\!+\!\e \frac{\bm w^TF_{\e}\bm y}{\bm w^T\bm y}{\bigg)}\!{\Bigg|}\!\leq\! \|\e\bm \zeta\|_1\right\}~~{\rm and}~~\left\{\la(M(\e))\!\setminus\!\breve{\la}(\e)\right\}\cap G_1\!=\!\varnothing.
\end{align}
\end{small}

We note that
\begin{align}\label{eq1242}
\bigg|\frac{\bm w^TF_{\e}\bm y}{\bm w^T\bm y}\bigg|\leq\max\{1,|\e|\}\cdot\underbrace{\left(\left|\frac{(\bm g^T\bm y_2)\cdot(\bm y_2^T \bm e)}{\bm y_1^T\bm y_2\D^2}-\frac{\bm y^T_2E\bm y_1 }{\bm y_1^T\bm y_2}\right|+\frac{(\bm y_2^T\bm e)^2}{2\D^2|\bm y_1^T\bm y_2|}\right)}_{\equiv\delta_0}.
\end{align}
Notice that $\widetilde{\nu}>0$. Denote by ${\la}_R(\e)$ the rightmost eigenvalue of $M(\e)$, and define
 $$
\varrho_0=\left\{
\begin{array}{ll}
1,               ~~~~~~~~~~~~~~~~~~~{\rm if}~\max\!\left\{\delta_0, \|W\|_1\right\}\!=\!0,\\
\specialrule{0em}{0.5ex}{0.5ex}
\frac{\widetilde{\nu}}{2\max\left\{\delta_0, \|W\|_1\right\}},~~ \rm else.
\end{array}\right.
$$
Next, we prove that
$$
{\la}_R(\e)=\breve{{\la}}(\e),~~{\rm for}~~ 0<|\e|<\widetilde{\rho}
=\min\{{\varrho}_0,\breve{\varrho},1\}
$$
by contradiction.
Suppose that ${\la}_R(\e)\neq \breve{{\la}}(\e)$ for $0<|\e|<\widetilde{\rho}$. It then follows from the second equality in \eqref{eq1451} that $\la_R(\e)\notin G_1$.
By the Ger$\rm \check{s}$gorin disc theorem \cite[Theorem 6.1.1]{Horn},
$$
\la_R(\e)\in \bigcup_{i =2}^{2n}G_i= \bigcup_{i=2}^{2n} {\Big\{}z:~{\big|}z-(t_{ii}+\e w_{ii}){\big|}\leq|\e\beta_i|+\sum_{j\neq i}|t_{ij}+\e w_{ij}|{\Big \}},
$$
where $\beta_i$ is the $i$-th element of $\bm \beta$ and $w_{ij}$ is the $(i,j)$-th element of $W$.
Thus, there is an integer $2\leq k\leq 2n$ such that $\la_R(\e)\in G_k$.

If $\max\left\{\delta_0, \|W\|_1\right\}=0$, then
$$
\mathfrak{Re}(t_{kk}+\e w_{kk})
 =\mathfrak{Re}(t_{kk})
<\la_*\overset{\eqref{eq1242}}{=}
 \la_*+\e \frac{\bm w^TF_{\e}\bm y}{\bm w^T\bm y},~~0<|\e|<\varrho_0=1.
$$
Otherwise, if $\max\left\{\delta_0, \|W\|_1\right\}\neq0$, as
$t_{kk}(\neq\la_*)$ is an eigenvalue of $M$, we have $\la_*-\mathfrak{Re}(t_{kk})\geq \widetilde{\nu}>0$. It follows that
 \begin{align*}
\mathfrak{Re}(t_{kk}\!+\!\e w_{kk})
\!\leq\!
 \la_*\!-\!\widetilde{\nu}+|\e~\mathfrak{Re}(w_{kk})|
&<\la_*\!-\! \widetilde{\nu}\!+\!\frac{\widetilde{\nu}~\|W\|_1}{2\max\left\{\delta_0, \|W\|_1\right\}}
\!\leq\!
 \la_*\!-\!\frac{\widetilde{\nu}}{2}\\
 &\leq
 \la_*\!+\!\e \frac{\bm w^T\!F_{\e}\bm y}{\bm w^T\!\bm y},{~~0\!<\!|\e|\!<\!\min\{\varrho_0,1\}}.
 \end{align*}
In conclusion, as $0<|\e|<\min\{\varrho_0,1\}$, we always  have $
\mathfrak{Re}(t_{kk}+\e w_{kk})<
 \la_*+\e \frac{\bm w^TF_{\e}\bm y}{\bm w^T\bm y}$. {That is,
the center $t_{kk}+\e w_{kk}$ of the circle  $G_k$ is to the left of the center $\la_{*}+\e\frac{\bm w^T\!F_{\e}\!\bm y}{\bm w^T\!\bm y}$ of the circle  $G_1$.}

 By Lemma \ref{lem20.59} (i), $\la_*(\e)\in G_k\cap\mathbb{R}$. Notice that $\la_*+\e \frac{\bm w^TF_{\e}\bm y}{\bm w^T\bm y}\in G_1\cap\mathbb{R}$ and $G_1 \cap G_k=\varnothing$ for $0<|\e|<\breve{\varrho}$. Thus, for $0<|\e|<\widetilde{\rho}
=\min\{{\varrho}_0,\breve{\varrho},1\}$, we know that  $G_1$ and $G_k$ satisfy the following two conditions:
\begin{equation*}
\left\{
\begin{array}{ll}
{\rm The~center~ of~the~circle}~ G_k~{\rm is~ to~ the~ left~ of~the~center~of~the~circle}~ G_1,\\
\specialrule{0em}{0.5ex}{0.5ex}
 G_1\cap \mathbb{R}\neq \varnothing,~G_k\cap \mathbb{R}\neq \varnothing~~{\rm and}~~
G_k\cap G_1=  \varnothing,
\end{array}\right.
\end{equation*}
as Fig. \ref{fig:3} depicts. Recall that $\breve{\la}(\e)\in G_1$, it is seen from Fig. \ref{fig:3} that $\la_R(\e)<\mathfrak{Re}(\breve{\la}(\e))$ as $0<|\e|<\widetilde{\rho}$, which contradicts to the fact that $\la_R(\e)$ is the rightmost eigenvalue of $M(\e)$.
Hence, we have $\la_R(\e)=\breve{\la}(\e)$ for $0<|\e|<\widetilde{\rho}$.
That is, $\breve{\la}(\e)$ is the rightmost eigenvalue $\la_R(\e)$ of $M(\e)$ for $0<|\e|<\widetilde{\rho}$. Moreover, by \cite[Proposition 3.4]{3}, $\la_R(\e)\in \mathbb{R}$.
Recall that $\la_*>0$ and the radius of $G_1$ is $\|\e\bm \zeta\|_1$, then
there is a $\varrho_1>0$ such that
$$\la_{*}+\e\frac{\bm w^TF_{\e}\bm y}{\bm w^T\!\bm y}>0~~{\rm and}~~
\la_{*}+\e\frac{\bm w^TF_{\e}\bm y}{\bm w^T\!\bm y}>\|\e\bm \zeta\|_1~~
{\rm for}~~|\e|\leq \varrho_1.
$$
Hence, $G_1\subseteq \{x: \mathfrak{Re}(x)>0\}$. It follows that $\la_R(\e)> 0$, and $\la_R(\e)=\la_*(\e)>0$ is the Lagrange multiplier of TRS \eqref{2} for
$\widetilde{\varrho}=\{\varrho_1,\widetilde{\rho}\}$.
\begin{figure}[H]
\begin{center}
\begin{tikzpicture}
\draw [ultra thick] (-2,0.8) circle [radius=1.78];
\draw [ultra thick] (2.5,0) circle [radius=1.08];
\draw [->](-4.5,0)--(-3,0)-- (-2,0)--(-1,0) -- (2,0) -- (4,0)--(5.5,0);
\node at (5.5,-0.3) {$\mathfrak{Re}$};
\node at (-2,1.06) {$t_{kk}+\e w_{kk}$};
\node at (-2.1,-0.3) {${\la}_{R}(\e)$};
\node at (-2,-1.5) {${G}_k$};
\node at (2.5,-0.35) {$\la_{*}\!+\!\e\frac{\bm w^T\!F_{\e}\!\bm y}{\bm w^T\!\bm y}$};
\node at (2.15,0.58) {$\breve{\la}(\e)$} ;
\node at (2.65,-1.5) {$G_1$};
\draw[fill](2.5,0) circle [radius=.8pt];
\draw[fill](-2.3,0) circle [radius=.8pt];
\draw[fill](-2,0.8) circle [radius=.8pt];
\draw[fill](2.15,0.25) circle [radius=.8pt];
\end{tikzpicture}
\end{center}\caption{}\label{fig:3}
\end{figure}

In summary, we have the first main theorem in this subsection:
\begin{theorem}\label{20.1200}
Suppose that $\la_*>0$, there is a scalar $\widetilde{\varrho}>0$ such that if $0<|\varepsilon|<\widetilde{\varrho}$, then
\begin{itemize}
\item[\rm(i)]
$\la_*(\e)>0$ is the rightmost eigenvalue of $M(\e)$ and
$\underset{\e\rightarrow 0}{\lim}{\la}_*(\e)=\la_*$,
\item[\rm(ii)]
 ${\la}_*(\varepsilon)$ is differentiable at $\varepsilon = 0$, and
\begin{eqnarray}\label{eq310}
{\la}'_*(0)=\frac{{\rm d}{\la}_*(\varepsilon)}{{\rm d}\varepsilon}{\bigg|_{\varepsilon=0}}=\frac{sign(\bm g^T\bm y_2)\|\bm y_1\|}{\bm y_1^T\bm y_2\D}\cdot\bm y_2^T\left(E\bm x_*+\bm e\right).
\end{eqnarray}
\end{itemize}
\end{theorem}

\begin{remark}
On one hand, as $\la_*$ is the rightmost eigenvalue of the matrix $M$, the condition number of the eigenvalue $\la_*$ can be defined as \cite{G.W.}
\begin{equation}
 cond(\la_{*})=\frac{1}{|\bm w^T\bm y|}=\frac{1}{2|\bm y_1^T\bm y_2|},
\end{equation}
where $\bm w$ is in \eqref{eq1023}.
On the other hand, in view of Theorem \ref{20.1200}, we give a {\tt new} definition on the condition number of $\la_*$. More precisely, note that if $\e$ is small enough, then
 $$
 \la_*(\e)=\la_*+\la'_*(0)\e+\mathcal{O}(\e^2).
 $$
 We have from \eqref{eq310} that
 \begin{align*}
 |\la_*'(0)|&={\bigg|}\frac{sign(\bm g^T\bm y_2)\|\bm y_1\|}{\D}\cdot\frac{\bm y_2^T\left(E\bm x_*+\bm e\right)}{\bm y_1^T\bm y_2}{\bigg |}\leq
  \frac{\|\bm y_1\|\|\bm y_2\|}{|\bm y_1^T\bm y_2|\D}\cdot(\|E\|\D+\|\bm e\|)\\
            &\leq \frac{\|\bm y_1\|\|\bm y_2\|}{\D|\bm y_1^T\bm y_2|}\cdot\max\{1,{\D}\}\cdot(\|E\|+\|\bm e\|).\nonumber
 \end{align*}
 Thus, if $\e$ small enough,
 \begin{align}\label{eq21.46}
 |\la_*-\la_*(\e)|\leq
 \frac{\|\bm y_1\|\|\bm y_2\|}{|\bm y_1^T\bm y_2|}\cdot\max\left\{1,\frac{1}{\D}\right\}\cdot\big(\|E\|+\|\bm e\|\big)|\e|+\mathcal{O}{\big(\e^2\big)}.
 \end{align}
Hence, we can define the condition number of $\la_*$ as follows
 \begin{equation}\label{eq3.9}
 s(\la_*)\!=\!\frac{\|\bm y_1\|\|\bm y_2\|}{|\bm y_1^T\!\bm y_2|}\cdot\max\left\{1,\frac{1}{\D}\right\}\!=\!\frac{\|A_*^{-1}\bm x_*\|}{\bm x_*^T\!A_*^{-1}\bm x_*}\cdot\max\left\{1,{\D}\right\},~~{ where}~~A_*= A+\la_*I,
 \end{equation}
and the last equality follows from \eqref{eq21.24} and \eqref{eq2039}.
Recall that neither $\bm y_1$ nor $\bm y_2$ are unit norm vectors, and $0<\|\bm y_1\|\|\bm y_2\|\leq\frac{1}{2}$. Therefore,
$$
s(\la_*)\leq cond(\la_{*})
$$
as $\D\geq 1$. Moreover, $\|\bm y_1\|\|\bm y_2\|$
can be much smaller than 1 in practice. Indeed,
$\|\bm y_1\|\|\bm y_2\|$ may be arbitrarily small when TRS is in ``nearly hard case" \cite{3}.
Consequently, $s(\la_*)$ can be much smaller than $cond(\la_{*})$; see Example 1 below.
\end{remark}

The following theorem establishes lower and upper bounds on $s(\la_*)$:
\begin{theorem}
Under the above notations, we have
$$
\max\left\{1,\frac{1}{\D} \right\}\leq s(\la_*)
   \leq \sqrt{\kappa(A_*)}\cdot \max\left\{1,\frac{1}{\D} \right\},
$$
where $A_*=A+\la_*I$ and  $\kappa(A_*)=\|A_*\|\|A_*^{-1}\|$.
\end{theorem}
\begin{proof}
From \eqref{eq16.12}, we have that
\begin{equation}\label{eq2039}
\bm y_1=A_*\bm y_2.
\end{equation}
Thus,
\begin{align*}
1&\leq
\frac{\|\bm y_1\|\|\bm y_2\|}{|\bm y_1^T\bm y_2|}=\frac{\|\bm y_1\|\|A_*^{-1}\bm y_1\|}{\bm y_1^TA_*^{-1}\bm y_1}\leq
\|A_*^{-\frac{1}{2}}\|\frac{\|\bm y_1\|\|A_*^{-\frac{1}{2}}\bm y_1\|}{\bm y_1^TA_*^{-1}\bm y_1}
      =\frac{\big\|A_*^{-\frac{1}{2}}\big\|\|\bm y_1\|}{\sqrt{\bm y_1^T A_*^{-1}\bm y_1}}\leq \!\sqrt{\kappa(A_*)},
\end{align*}
which completes the proof.
\end{proof}

Second, we focus on the continuity and differentiability of $\bm x_*(\e)$ at $\e=0$.
The following equations are needed \cite[p. 52]{Y.S}:
\begin{subequations}
\begin{align}
&(B_1-z I)^{-1}-({B}_2-z I)^{-1} = (B_1-z I)^{-1}(B_2-{B_1})({B}_2-zI)^{-1},\label{8.2777}\\
&(B-z_1 I)^{-1}- (B-z_2 I)^{-1}= (z_2-z_1)(B-z_1 I)^{-1}(B-z_2 I)^{-1}\label{8.25}
\end{align}
\end{subequations}
where $B,B_1,B_2\in \mathbb{C}^{n\times n}$ and $z,z_1,z_2\in \mathbb{C}$. We have the following lemma.
\begin{lemma}\label{lem1531}
Suppose that TRS \eqref{1} is in easy case,
there is a scalar ${\varrho}>0$ such that if
  $0<|\e|<{\varrho}$, then TRS \eqref{2} is also in easy case and
\begin{equation}\label{eq16.24}
\begin{aligned}
&\left(I+(\la_*(\e)-\la_*)(A(\e)+\la_* I)^{-1}\right)\cdot(\bm x_*(\e)-\bm x_*)\\
=&
-(\la_*(\e)-\la_*) ({A}(\e)+{\la}_{*}I)^{-1}{\bm x_*}-\e\cdot({A(\e)}+{\la}_{*}I)^{-1}(E\bm x_*+\bm e).
\end{aligned}
\end{equation}
\end{lemma}
\begin{proof}
As TRS \eqref{1} is in \emph{easy case}, by Theorem \ref{thm2119}, there is a scalar $\varsigma_1>0$, such that  TRS \eqref{2} is also in \emph{easy case},
 i.e., $A(\e)+\la_*(\e)I\succ\bm O$
 as $0<|\e|<\varsigma_1$.
 From Theorem \ref{7.59}, we have that
 $$
 {\bm x}_{*}(\e)
=-\big({A}(\e)+{\la}_{*}(\e)I\big)^{-1}{\bm g(\e)}~~{\rm and}~~\bm x_{*}=-(A+\la_{*}I)^{-1}\bm g.
 $$
If $0<|\e|<\delta_1=\frac{1}{\|(A+\la_*I)^{-1}E\|}$, then
${A}(\varepsilon)+\la_{*}I=(A+\la_*I)\left(I+\e (A+\la_*I)^{-1}E\right)$
is nonsingular. Thus, if $0<|\e|<\varrho=\min\{\varsigma_1, \delta_1\}$, there holds
\begin{align}
&~~~~{\bm x}_{*}(\e)-\bm x_{*}
=(A+\la_{*}I)^{-1}\bm g-\big({A}(\e)+{\la}_{*}(\e)I\big)^{-1}{\bm g(\e)}\nonumber\\
 &=\big({A}(\e)+{\la_*}I\big)^{-1}{\bm g}-\big({A}(\e)+{\la}_{*}(\e)I\big)^{-1}{\bm g(\e)}-\big[\big({A}(\e)+{\la}_{*}I\big)^{-1}{\bm g}-(A+\la_{*}I)^{-1}\bm g\big]\nonumber\\
 &\overset{\eqref{8.2777}}{=} \big({A}(\e)+{\la_*}I\big)^{-1}{\bm g(\e)}-\e\cdot\big({A}(\e)+{\la_*}I\big)^{-1}\bm e-\big({A}(\e)+{\la}_{*}(\e)I\big)^{-1}{\bm g(\e)}\nonumber\\
 &~~~~+\e\cdot({A}(\e)+{\la}_{*}I)^{-1}E(A+\la_{*}I)^{-1}{\bm g}\nonumber\\
  &\overset{\eqref{8.25}}{=}(\la_*(\e)-\la_*) \cdot \big({A}(\e)+{\la}_{*}I\big)^{-1}\big({A}(\e)+{\la_*(\e)}I\big)^{-1}{\bm g(\e)}-\e\cdot(A(\e)+\la_*I)^{-1}\bm e\nonumber\\
  &~~~~-\e\cdot({A(\e)}+{\la}_{*}I)^{-1}E\bm x_*\nonumber\\
  &=-(\la_*(\e)-\la_*) \big({A}(\e)+{\la}_{*}I\big)^{-1}{\bm x_*(\e)}-\e\cdot({A(\e)}+{\la}_{*}I)^{-1}(E\bm x_*+\bm e)\nonumber\\
  &=-(\la_*(\e)-\la_*) \big({A}(\e)+{\la}_{*}I\big)^{-1}{\bm x_*}-(\la_*(\e)-\la_*)\big({A}(\e)+{\la}_{*}I\big)^{-1}(\bm x_*(\e)-{\bm x_*})\nonumber\\
  &~~~~-\e\cdot({A(\e)}+{\la}_{*}I)^{-1}(E\bm x_*+\bm e),\nonumber
\end{align}
which completes the proof.
\end{proof}

With the help of Lemma \ref{lem1531}, we have the second main theorem in this subsection, which shows the continuity and differentiability of $\bm x_*(\e)$ at $\e=0$.
\begin{theorem}\label{coro1532}
Suppose that TRS \eqref{1} is in easy case and $\la_{*}>0$,
 there exists a scalar $\underline{\varrho}>0$, such that if $0<|\varepsilon|<\underline{\varrho}$, then
\begin{itemize}
\item[(i)]
$\|\bm x_*(\e)\|=\|\bm x_*\|=\D$
 and
$\underset{\e\rightarrow 0}{\lim} \bm x_*{(\e)} = \bm x_*$,
\item[(ii)]
 ${\bm x}_*(\varepsilon)$ is differentiable at $\varepsilon = 0$, and
\begin{align}\label{eq20270}
\!\!\!\!\bm x'_*(0)=\frac{{\rm d}\bm x_*(\e)}{{\rm d}\e}{\bigg|_{\e=0}}\!=\!
\left(-A_*^{-1}\!+\!\frac{\bm y_2\bm y_2^T}{\bm y_1^T\bm y_2}\right)(E\bm x^*\!+\!\bm e),~~where~~A_*\!=\!A\!+\!\la_*I.
\end{align}
\end{itemize}
\end{theorem}
\begin{proof}
Notice that $\|\bm x_*\|=\D$ as $\la_{*}>0$.
We obtain from Theorem \ref{20.1200} that $\la_*(\e)>0$ if $0<|\e|<\widetilde{\varrho}$, and $\|\bm x_*(\e)\|=\D$ from \eqref{eq2100}.
Let $0<|\e|<\underline{\varrho}=\min\{\varrho,\widetilde{\varrho}\}$, it follows from \eqref{eq16.24} and Theorem \ref{20.1200} (i) that $\underset{\e\rightarrow 0}{\lim} \bm x_*{(\e)} = \bm x_*$, and
\begin{align}
&\lim_{\e\rightarrow 0}\Big(I+(\la_*(\e)-\la_*)\big(A(\e)+\la^* I\big)^{-1}\Big)\cdot \lim_{\e\rightarrow 0}\frac{\bm x_*(\e)-\bm x_*}{\e}
\nonumber\\
=&-\lim_{\e\rightarrow 0}({A(\e)}+{\la}_{*}I)^{-1}(E\bm x_*\!+\!\bm e)\!-\!\lim_{\e\rightarrow 0}\frac{\la_*(\e)\!-\!\la_*}{\e} \cdot\lim_{\e\rightarrow 0}\big({A}(\e)\!+\!{\la}_{*}I\big)^{-1}\!{\bm x_*}\label{eq1659}\\
=&-A_*^{-1}(E\bm x^*+\bm e)-\la'_*(0)\cdot A_*^{-1}\bm x_*\nonumber\\
=&-A_*^{-1}(E\bm x^*+\bm e)-{\bigg (}\frac{sign(\bm g^T\bm y_2)\|\bm y_1\|}{\D}\cdot\frac{\bm y_2^T\left(E\bm x_*+\bm e\right)}{\bm y_1^T\bm y_2}{\bigg )} A_*^{-1}\bm x_*\nonumber\\
\overset{\eqref{eq21.24}}{=}&\left(-A_*^{-1}+\frac{A_*^{-1}\bm y_1\bm y_2^T}{\bm y_1^T\bm y_2}\right)(E\bm x^*+\bm e)
\overset{\eqref{eq2039}}{=}\left(-A_*^{-1}+\frac{\bm y_2\bm y_2^T}{\bm y_1^T\bm y_2}\right)(E\bm x^*+\bm e).\nonumber
\end{align}
By Theorem \ref{20.1200} (i), $\underset{\e\rightarrow 0}{\lim}\left(I+(\la_*(\e)-\la_*)\big(A(\e)+\la^* I\big)^{-1}\right)=I$, and thus
$$
\bm x'(0)=\lim_{\e\rightarrow 0}\frac{\bm x_*(\e)-\bm x_*}{\e}=\left(-A_*^{-1}+\frac{\bm y_2\bm y_2^T}{\bm y_1^T\bm y_2}\right)(E\bm x^*+\bm e).
$$
\end{proof}
\begin{remark}
We are in a position to define the condition number of $\bm x_*$.
If $\e$ is small enough, then
$$
\bm x_*(\e)=\bm x_*+\bm x_*'(0) \e+\mathcal{O}(\e^2).
$$
It is seen from \eqref{eq20270} and $\D=\|\bm x_*\|=\|A^{-1}_*\bm g\|\geq \frac{\|\bm g\|}{\|A_*\|}$ that
\begin{align}
\frac{\|\bm x_*(\e)-\bm x_*\|}{\D}\leq& \frac{\|\bm x_*'(0)\|}{\D} |\e|+\mathcal{O}\big(\e^2\big)
\leq\left\|A_*^{-1}-\frac{\bm y_2\bm y_2^T}{\bm y_1^T\bm y_2}\right\|\cdot\left(\|E\|+\frac{\|\bm e\|}{\D}\right)|\e|+\mathcal{O}\big(\e^2\big)\nonumber\\
 \leq&\|A_{*}\|\left\|A_*^{-1}-\frac{\bm y_2\bm y_2^T}{\bm y_1^T\bm y_2}\right\|\cdot\left(\frac{\|E\|}{\|A_*\|}+\frac{\|\bm e\|}{\|\bm g\|}\right)|\e|+\mathcal{O}\big(\e^2\big).\label{eq1531}
\end{align}
As a result, we define
\begin{equation}\label{eq3.16}
s(\bm x_*)=\|A_{*}\|\left\|A_*^{-1}-\frac{\bm y_2\bm y_2^T}{\bm y_1^T\bm y_2}\right\|{=\|A_{*}\|\left\|A_*^{-1}- \frac{A_*^{-1}\bm x_*\bm x_*^TA_*^{-1}}{\bm x_*^TA_*^{-1}\bm x_*}\right\|}
\end{equation}
as the condition number of $\bm x_*$, where the second equality follows from
\eqref{eq21.24} and \eqref{eq2039}.
\end{remark}
\begin{lemma}\cite[Theorem 3.1.3]{GenInv}\label{lemm}
For $W\in \mathbb{R}^{n\times n}$, $\bm c\in\mathbb{R}^n$ and $\bm d\in \mathbb{R}^n$, let $\bm f = W^\dagger\bm c$ and
$\bm h = (W^\dagger)^T\bm d$.
 If $\bm c\in \mathcal{R}(W)$, $\bm d\in\mathcal{R}(W^T)$ and
 $1+\bm c^TW^\dag\bm d=0$, then
 $$
 (W+\bm c\bm d^T)^\dag= W^\dag-\frac{\bm f\bm f^T W^\dag}{\|\bm f\|^2}-\frac{W^\dag\bm h\bm h^T}{\|\bm h\|^2}+\frac{\bm f\bm f^T W^\dag\bm h\bm h^T}{(\|\bm h\|\|\bm f\|)^2}.
 $$
\end{lemma}
The following theorem establishes lower and upper bounds on $s(\bm x_*)$:
\begin{theorem}\label{thm111}
Under the above notations, we have that
\begin{itemize}
\item[\rm (i)] Let $P= I-\frac{\bm x_*\bm x_*^T}{\D^2}$, then
\begin{align}A_*^{-1}- \frac{A_*^{-1}\bm x_*\bm x_*^TA_*^{-1}}{\bm x_*^TA_*^{-1}\bm x_*}
=(PA_*P)^\dag,\nonumber
\end{align}
and
$s(\bm x_*)=\|A_{*}\|\left\|(PA_*P)^\dag\right\|$;
\item[\rm(ii)]
0 is a simple eigenvalue of $PA_*P$;
\item[\rm (iii)]
$1\leq \frac{\alpha_1+\la_*}{\alpha_{n-1}+\la_*}\leq s(\bm x_*)\leq \frac{\alpha_1+\la_*}{\alpha_{n}+\la_*}= \kappa(A_*)$.
\end{itemize}
\end{theorem}
\begin{proof}
(i) Let
$W=A_*^{-1}, ~~\bm c=\frac{A_*^{-1}\bm x_*}{\sqrt{\bm x_*^TA_*^{-1}\bm x_*}}~~ {\rm and}~~ \bm d=\frac{-A_*^{-1}\bm x_*}{\sqrt{\bm x_*^TA_*^{-1}\bm x_*}}.$
Thus,
$$
\bm c\in \mathcal{R}(W),~~\bm d\in\mathcal{R}(W^T),~~{\rm and}~~
1+\bm c^TW^\dag\bm d =1- \frac{\bm x_*^TA_*^{-1}A_*A_*^{-1}\bm x_*}{\bm x_*^TA_*^{-1}\bm x_*}=0.
$$
By Lemma \ref{lemm},
\begin{align*}&\left(A_*^{-1}- \frac{A_*^{-1}\bm x_*\bm x_*^TA_*^{-1}}{\bm x_*^TA_*^{-1}\bm x_*}\right)^\dag
\!\!=\!A_*-\frac{\bm x_*\bm x_*^TA_*+A_*\bm x_*\bm x_*^T}{\|\bm x_*\|^2} +\frac{(\bm x_*^TA_*\bm x_*)\bm x_*\bm x_*^T}{\|\bm x_*\|^4}
=PA_*P,
\end{align*}
which gives (i).

(ii) Notice that (0,~$\bm x_*$) is an eigenpair of $PA_*P$ which is symmetric. If 0 is a multiple eigenvalue, then there is a vector $\bm a\neq \bm 0$ with ${\bm x}^T_*\bm a=0$, such that
$$
\bm a^TPA_*P\bm a=\bm a^T\left(I-\frac{\bm x_*\bm x_*^T}{\D^2}\right)A_*\left(I-\frac{\bm x_*\bm x_*^T}{\D^2}\right)\bm a = 0.
$$
It follows that $\bm a^TA_*\bm a=0$, which contradicts to the fact that $A_*\succ\bm O$.
\\
(iii) Let $\mathcal{A}=A_*^{-1}-\frac{\bm y_2\bm y_2^T}{\bm y_1^T\bm y_2}$, where $\bm y_1$ and $\bm y_2$ are defined in Theorem \ref{16.0200}.
 By (i), $\mathcal{A}\succcurlyeq\bm O$ and $\|\mathcal{A}\|$ is the largest eigenvalue of $\mathcal{A}$.
As $A_*^{-1}\succ\bm O$, $\frac{\bm y_2\bm y_2^T}{\bm y_1^T\bm y_2}=\frac{\bm y_2\bm y_2^T}{\bm y_2^TA_*\bm y_2}\succcurlyeq\bm O$ and $rank(\bm y_2\bm y_2^T)=1$,
it follows from \cite[Corollary 4.3.9]{Horn} that
$$ (\alpha_{n-1}+\la_*)^{-1}\leq\|\mathcal{A}\|\leq\|{A}_*^{-1}\|,$$
from which we get (iii).
\end{proof}
\begin{example}
Counter-intuitively, we illustrate that $\la_*$ and $\bm x_{*}$ may be well-conditioned even if TRS is in ``nearly hard case'' or $A_*$ is ill-conditioned in this example. Consider
$$
A = \begin{pmatrix}
1&\\
  &-1+\gamma
\end{pmatrix},~\bm g = \begin{pmatrix}
1\\
\gamma
\end{pmatrix}~~{and}~~\D = \frac{\sqrt{5}}{2},
$$
where $0<\gamma\ll 1$.
We have that
$
\la_*=1,~\bm x_*= (
-0.5~
-\!1
)^T~~{and}~~A_*=diag(
2, \gamma).
$
Moreover,
$$
\left\{
\begin{array}{ll}
s(\la_*)  &=\frac{\|A_*^{-1}\bm x_*\|}{\bm x_*^TA_*^{-1}\bm x_*}\max\left\{1,{\D}\right\}=\frac{\sqrt{5}\left\|\begin{pmatrix}
\frac{1}{4} &\gamma^{-1}
\end{pmatrix}^T\right\|}{2\big(\frac{1}{8}+\gamma^{-1}\big)}
=\frac{\sqrt{5}\sqrt{16+\gamma^2}}{8+\gamma}<\frac{6}{5},\\
\specialrule{0em}{1ex}{1ex}
s(\bm x_*)&=\|A_*\|\|(PA_*P)^\dag\|=\frac{2}{8+\gamma}\cdot\left\|\begin{pmatrix}
4 &-2\\
-2 &1
\end{pmatrix}
\right\|=\frac{10}{8+\gamma}<\frac{5}{4}.
\end{array}\right.
$$
As $\gamma$ is close to 0, $A_*$ is a nearly singular matrix, and TRS is in ``nearly hard case'' \cite{3,M1}.

Let $E = diag(-1,1)$, $\bm e = (
1~~ 2 )^T$, and $0<\e<\gamma$. Then
$$
A(\e) =A+\e E= \begin{pmatrix}1-\e &\\ &-1+\gamma+\e
\end{pmatrix}~~{and}~~\bm g(\e)=\bm g+\e\bm e = \begin{pmatrix}1+\e\\ \gamma+2\e\end{pmatrix}.
$$
Denote
$$
\wp(\la)\equiv \left\|\left(A(\e)+\la I\right)^{-1}\bm g(\e)\right\|^2=
\frac{(1+\e)^2}{(\la+1-\e)^2}+\frac{(\gamma+2\e)^2}{(\la-1+\gamma+\e)^2},
$$
we see that
$$
\wp(1+\e)\geq \D^2\geq \wp(1+\gamma\e+{\e}).
$$
Note that $A(\e)\!+\!\la_*(\e)I\succ\bm O$, $\wp\left(\la_*(\e)\right)\!=\!\D^2$ and
 $\wp(\la)$ is a continuous and strictly monotonically decreasing function on $\la\!\in\![1\!+\!\e,1\!+\!{\e}\!+\!\gamma\e]$. Consequently, $\la_*(\e)\!\in\![1\!+\!\e,1\!+\!{\e}\!+\!\gamma\e]$, i.e.,
 $$
 \la_*(\e)=1+ p\e=\la_*+p\e~~{with}~~p\in[1,1+\gamma].
 $$
 It follows that
 $$
 \frac{\|\bm x_*(\e)\!-\!\bm x_*\|}{\D}\!=\!\frac{\sqrt{5}} {2}\!\cdot\!\left\|\left(A(\e)\!+\!\la_*(\e) I\right)^{-1}\!\!\bm g(\e)\!-\!\begin{pmatrix}
\frac{1}{2} \\ 1
\end{pmatrix}\right\|\! =\!\widetilde{p} \e, ~~{with}~~
 \widetilde{p}\!\in\!\!\left[\frac{\sqrt{5}(2\!-\!\gamma)}{4(2\!+\!\gamma)} ,\frac{5}{4}\right].
 $$
It is observed from Table \ref{T3.2} that our bounds \eqref{eq21.46} and \eqref{eq1531} are very sharp. Moreover, $s(\la_*)$ and $s(\bm x_{*})$ may not be large even if TRS is in ``nearly hard case'' or $\kappa(A_*)$ is large.
In summary, $\la_*$ and $\bm x_{*}$ could be well-conditioned even if TRS is in ``nearly hard case'' or $A_*$ is ill-conditioned, and $\kappa(A_*)$ is {\tt unsuitable} to depict the conditioning of TRS.
 \begin{table}[H]
\begin{center}\caption{Sharpness of \eqref{eq21.46} and \eqref{eq1531}}\label{T3.2}
\begin{tabular}{c|c}
\hline
\tabincell{c}{$|\la_*(\e)-\la_*|=p\e$,\\ with $p\!\in\![1,1+\gamma]$}&\tabincell{c}{ Bound in \eqref{eq21.46}\!~:~$\frac{\sqrt{16+\gamma^2}}{8+\gamma}(5+\sqrt{5})\e+\mathcal{O}(\e^2)
\leq 3\e $} \\
   \hline
\tabincell{c}{$\frac{\|\bm x_*(\e)-\bm x_*\|}{\D}=\widetilde{p} \e$, \\ with
 $\widetilde{p}\in\big[\frac{\sqrt{5}(2-\gamma)}{4(2+\gamma)},\frac{5}{4}\big]$} &\tabincell{c}{Bound in \eqref{eq1531}\!~:~ $\frac{10}{8+\gamma}\Big(
\frac{1}{2}\!+\!\sqrt{\frac{5}{1+\gamma^2}}\Big)\e+\mathcal{O}(\e^2)
\leq \frac{7}{2}\e$ }\\
   \hline
\end{tabular}
\end{center}
\end{table}
Finally, we notice that
\begin{align*}
cond(\la_*)\!=&\frac{1}{2|\bm y_1^T\!\bm y_2|}\overset{ \eqref{eq3.9}}{=}\frac{s(\la_*)}{2\|\bm y_1\|\|\bm y_2\|}
\!>\!\frac{\|A_*^{-1}\bm y_1\|}{2\|\bm y_1\|}\overset{\eqref{eq21.24}}{=}\frac{\|A_*^{-1}\bm x_*\|}{2\|\bm x_*\|}
=\frac{1}{\sqrt{5}}\left\|\begin{bmatrix}
0.25\\ \gamma^{-1}
\end{bmatrix}\right\|\!>\!\frac{\gamma^{-1}}{\sqrt{5}},
\end{align*}
where the first inequality is from $s(\la_*)\geq1>\|\bm y_2\|^2\overset{\eqref{eq2039}}{=}\|\bm y_2\|{\|A_*^{-1}\bm y_1\|}$.
That is, $cond(\la_*)\gg s(\la_*)$ and $cond(\la_*)\gg 1$ as $\gamma\ll 1$.
In conclusion, although $\la_*$ is an eigenvalue of $M$, the eigenvalue condition number $cond(\la_*)$
is inappropriate to access the conditioning of $\la_*$ in TRS \eqref{1}.
\end{example}

\subsection{The case of $\bm{\la_{*}=0}$ and $\bm{\| x_{*}\|=\D}$}\label{sec14}
In this case, we point out that $\la_*(\e)$ and $\bm x_*(\e)$ may not be differentiable at $\e=0$.
For instance, consider
$$
A=\begin{pmatrix}
2 & \\
   & 1
\end{pmatrix},~
\bm g=\begin{pmatrix}
1 \\ 1
\end{pmatrix}
,~
\D=\frac{\sqrt{5}}{2},~E=\begin{pmatrix}
1 & \\
   & 0
\end{pmatrix}~~{\rm and}~~\bm e=\begin{pmatrix}
0\\
    -1
\end{pmatrix}.
$$
We have $\bm x_{*}=A^{-1}\bm g=\begin{pmatrix}
-\frac{1}{2} &-1
\end{pmatrix}^T$, $\la_*=0$, and
$$
A(\e)=\begin{pmatrix}
2+\e & \\
   & 1
\end{pmatrix},\quad
\bm g(\e)=\begin{pmatrix}
1 \\ 1-\e
\end{pmatrix},
$$
where $0<|\e|< 1$.
We have
$$
\left\{
\begin{array}{ll}
\la_*(\e)=0~~~{\rm and}~~~\bm x_*(\e)=(-\frac{1}{2+\e}~~-\!\!1\!+\!\e)^T~~                    &{\rm if}~~ \e>0,\\
\specialrule{0em}{0.2ex}{0.2ex}
\la_*(\e)=\e~~~{\rm and}~~~\bm x_*(\e)=(-\frac{1}{2}~~-\!1)^T~~ &{\rm if}~~ \e<0.
\end{array}\right.
$$
Thus,
$$\underset{\e\rightarrow 0}{\lim}\la_*(\e)=\la_*=\la_*(0)~~~{\rm and}~~~\underset{\e\rightarrow 0}{\lim} \bm x_*{(\e)} = \bm x_*=\bm x_*(0).$$
 However,
 $$\left\{
\begin{array}{ll}
-1=\underset{\e\rightarrow 0^-}{\lim}\frac{\la_*(\e)-\la_*}{\e}\neq\underset{\e\rightarrow 0^+}{\lim}\frac{\la_*(\e)-\la_*}{\e}=0,
\\
\specialrule{0em}{0.8ex}{0.8ex}
\bm 0=\underset{\e\rightarrow 0^-}{\lim}\frac{\bm x_*(\e)-\bm x_*}{\e}\neq\underset{\e\rightarrow 0^+}{\lim}\frac{\bm x_*(\e)-\bm x_*}{\e}=\left(
\frac{1}{4}~~ 1\right)^T.
\end{array}\right.
$$
Hence, $\la_*(\e)$ and $\bm x_*(\e)$ are not differentiable at $\e=0$.
The following theorem gives a sufficient and necessary condition for $\la_*(\e)$ and $\bm x_*(\e)$ being {\it non-differentiable} at $\e=0$.
More precisely, it shows that if $ \la_{*}=0,~\|\bm x_{*}\|=\D$, then $\la_*(\e)$ and $\bm x_*(\e)$ are \emph{non-differentiable} at $\e=0$ {\it if and only if}
$\bm x_*^TA^{-1}(\bm e+E\bm x_*)\neq 0$.

\begin{theorem}\label{Thm3.13}
Suppose that \emph{TRS \eqref{1}} is in easy case, and $\|\bm x_{*}\|=\D$, $\la_{*}=0$. We have that
  \begin{itemize}
 \item[\rm(i)] $\la_*(\e)$ and $\bm x_*(\e)$ are continuous at $\e = 0$ and $\underset{\e\rightarrow 0}{\lim}\la_*(\e)=\la_*$, $\underset{\e\rightarrow 0}{\lim} \bm x_*{(\e)} = \bm x_*$,
 \item[\rm(ii)] if $\bm x_*^TA^{-1}(\bm e+E\bm x_*)\neq 0$, then both $\la_*(\e)$ and $\bm x_*(\e)$ are {non-differentiable} at $\e=0$,
 \item[\rm(iii)]
      if $\bm x_*^TA^{-1}(\bm e+E\bm x_*)= 0$, then both $\la_*(\e)$ and $\bm x_*(\e)$ are differentiable at $\e=0$, moreover, $\la'(0)=0$ and $\bm x'_*(0)=-A^{-1}(E\bm x^*+\bm e)$.
 \end{itemize}
\end{theorem}
\begin{proof}
We first prove (i) and (ii). Notice that $A=A+\la_{*}I\succ\bm O$ as $\la_*=0$, and if $|\e|\|E\|<\alpha_n$, we have $A(\e)\succ\bm O$.
As $\varepsilon\rightarrow 0$, without loss of generality,
we assume that $|\e|<\frac{\alpha_n}{\|E\|}$ such that $A(\e)$ is positive definite.
Denote by $\widetilde{\bm x}= -A(\e)^{-1}\bm g(\e)$, then it follows from the classical perturbation theorem on linear system that \cite[eq. (1.75)]{2}
 $$
 \frac{{\rm d}\widetilde{\bm x}}{{\rm d} \e}\bigg|_{\e=0}=-A^{-1}(E\bm x^*+\bm e).
 $$
Denote by $\phi(\e)=\|\widetilde{\bm x}\|^2$, we have that $\phi(0)=\|\bm x_*\|^2=\D^2$.
Thus,
\begin{equation}\label{eq22.22}
\|\widetilde{\bm x}\|^2=\phi(0)+\phi'(0)\e+\frac{\phi''(0)}{2}\e^2+\mathcal{O}(\e^3)
=\D^2+\phi'(0)\e+\frac{\phi''(0)}{2}\e^2+\mathcal{O}(\e^3),
\end{equation}
where
$$
\phi'(0)=\frac{{\rm d}\|\widetilde{\bm x}\|^2}{{\rm d} \e}\bigg|_{\e=0}
=\frac{{\rm d}(\widetilde{\bm x}^T\widetilde{\bm x})}{{\rm d} \e}\bigg|_{\e=0}
=-2\bm x_*^TA^{-1}(\bm e+E\bm x_*),
$$
and $\phi''(0)$ is uniformly bounded.

If $\bm x_*^TA^{-1}(\bm e+E\bm x_*)>0$, there is a scalar $\delta^+>0$ such that for any $\e\in (0,~\delta^+]$, $\|\widetilde{\bm x}\|< \D$; refer to \eqref{eq22.22}. By Theorem \ref{7.59}, we have $\widetilde{\bm x}=\bm x_*(\e)$ and $\la_*(\e)=\la_*=0$. Therefore,
$$
\lim_{\e\rightarrow 0^+}\bm x_*(\e)=\bm x_*,~~\lim_{\e\rightarrow 0^+}\la_*(\e)=\la_*~~{\rm and}~~ \lim_{\e\rightarrow 0^+}\frac{\la_*(\e)-\la_*}{\e}=0.
$$
On the other hand, there is a ${\delta}^-<0$ such that for any $\e\in [{\delta}^-,0)$, $\|\widetilde{\bm x}\|>\D$. So we have  $\la_*(\e)>0$ and thus $\|\bm x_*(\e)\|=\D$. By Theorem \ref{16.0200}, $\la_*(\e)$ is the rightmost eigenvalue of $M(\e)$ for all $\e\in [{\delta}^-,0)$.
Recall that $\la_*$ is the rightmost eigenvalue of $M$ if $\|\bm x_*\|=\D$. From \eqref{eq16.24} and Theorem  \ref{20.1200},
we obtain
\begin{small}
 \begin{equation*}
\lim_{\e\rightarrow 0^-}\bm x_*(\e)\!=\!\bm x_*,~\lim_{\e\rightarrow 0^-}\la_*(\e)\!=\!\la_*~
~{\rm and}~
\lim_{\e\rightarrow 0^-}\frac{\la_*(\e)-\la_*}{\e}\!=\!\frac{sign(\bm g^T\bm y_2)\|\bm y_1\|}{\bm y_1^T\bm y_2\D}\cdot\bm y_2^T\left(E\bm x_*\!+\!\bm e\right).
\end{equation*}
\end{small}
As a result, we have $\underset{\e\rightarrow 0}{\lim}\bm x_*(\e)\!=\!\bm x_*$ and $\underset{\e\rightarrow 0}{\lim}\la_*(\e)\!=\!\la_*$.

Moreover, it follows from  \eqref{eq2039} that $\bm y_1^T\bm y_2>0$, and
\begin{align*}
&\bm x_*^TA^{-1}(\bm e\!+\!E\bm x_*)\!\overset{\eqref{eq21.24}}{=}\!\frac{-sign(\bm g^T\!\bm y_2)\D}{\|\bm y_1\|}\cdot\bm y^T_2
(\bm e\!+\!E\bm x_*)\!=\!-\frac{\bm y_1^T\bm y_2\D^2}{\|\bm y_1\|^2}\cdot\left[\lim_{\e\rightarrow 0^-}\!\frac{\la_*(\e)\!-\!\la_*}{\e}\right].
\end{align*}
However,
$$
0=\lim_{\e\rightarrow 0^+}\frac{\la_*(\e)-\la_*}{\e}\neq \lim_{\e\rightarrow 0^-}\frac{\la_*(\e)-\la_*}{\e}=-\frac{\|\bm y_1\|^2\bm x_*^TA^{-1}(\bm e+E\bm x_*)}{\bm y_1^T\bm y_2\D^2}<0.
$$
Thus, $\la_*(\e)$ is non-differentiable at $\e = 0$. Similarly, we have from \eqref{eq1659} that $\bm x_*(\e)$ is non-differentiable at $\e = 0$.
By using the same trick, when $\bm x_*^TA^{-1}(\bm e+E\bm x_*)<0$, we can prove that $\bm x_*(\e)$ and $\la_*(\e)$ are also non-differentiable at $\e=0$.


Next, we will prove (iii). Suppose that $\bm x_*^TA^{-1}(\bm e+E\bm x_*)=0$.
On one hand, if $\|\widetilde{\bm x}\|=\|A(\e)^{-1}\bm g(\e)\|<\D$, we have from Theorem \ref{7.59} that $\la_*(\e)=\la_*=0$ and $\bm x_*(\e)=-A(\e)^{-1}\bm g(\e)$.
Then, it follows that $\underset{\e\rightarrow 0}{\lim}\bm x_*(\e)=\bm x_*$, $\la_*'(0)=0$ and $\bm x'_*(0)=-A^{-1}(E\bm x^*+\bm e)$ \cite[eq. (1.75)]{2}.
On the other hand, if  $\|\widetilde{\bm x}\|=\|A(\e)^{-1}\bm g(\e)\|\geq\D$, we have from \eqref{eq22.22} that there exists a constant $d$ independent of $\e$, such that
\begin{equation}\label{eq2018}
\D^2\leq\|\widetilde{\bm x}\|^2\leq\D^2 + d\e^2~~{\rm with}~~0\leq d=\mathcal{O}\big(|\phi''(0)|\big)<\infty.
\end{equation}

Let
$$
p_{\e}(\la)\!=\!\big\|(A(\e)\!+\!\la I)^{-1}\!\bm g(\e)\big\|^2\!=\!\sum_{i=1}^{n}\frac{\big(\widetilde{\bm u}_i^T\!\bm g(\e)\big)^2}{(\la\!+\!\widetilde{\alpha}_i)^2},~~{\rm with}~~\la\!\in\![0, +\infty)~~{\rm and}~~0\!<\!|\e|\!<\!\frac{\alpha_n}{\|E\|}.
$$
Thus, $p_{\e}(\la)$ is a continuous and strictly monotonically decreasing function on $\la$. By Taylor's theorem \cite[Section 5.4]{Mathanal}, there is
a scalar $\vartheta_{\e}\in [0,\la]$, such that
\begin{align*}
p_{\e}(\la)=&p_{\e}(0)+p_{\e}'(0)\la+\frac{1}{2}p_{\e}''(\vartheta_{\e})\la^2\\
=&\|\widetilde{\bm x}\|^2-2\la\cdot\sum_{i=1}^{n}\frac{\big(\widetilde{\bm u}_i^T\bm g(\e)\big)^2}{\widetilde{\alpha}_i^3}+3\la^2\cdot\sum_{i=1}^{n}
\frac{\big(\widetilde{\bm u}_i^T\bm g(\e)\big)^2}{(\vartheta_{\e}+\widetilde{\alpha}_i)^4}\\
=&\|\widetilde{\bm x}\|^2-2\la\cdot\bm g(\e)^TA(\e)^{-3}\bm g(\e)+3\la^2\cdot\bm g(\e)^T(A(\e)+\vartheta_{\e} I)^{-4}\bm g(\e).
\end{align*}
As
 $$\xi\e^2=\frac{d}{\bm g(\e)^TA(\e)^{-3}\bm g(\e)}\cdot\e^2\geq 0,~~{\rm if}~~0<|\e|<\frac{\alpha_n}{\|E\|},$$
we have
\begin{align}\label{eq1416}
p_{\e}\left(\xi\e^2\right)=\|\widetilde{\bm x}\|^2-2d\e^2+3d^2\e^4\cdot\frac{\bm g(\e)^T(A(\e)+\vartheta_{\e} I)^{-4}\bm g(\e)}{\big(\bm g(\e)^TA(\e)^{-3}\bm g(\e)\big)^2}.
\end{align}
Note that
$$
\bm g(\e)^T(A(\e)+\vartheta_{\e} I)^{-4}\bm g(\e)\leq \bm g(\e)^TA(\e)^{-4}\bm g(\e),
$$
and
$$
\underset{\e\rightarrow 0}{\lim}\bm g(\e)^TA(\e)^{-i}\bm g(\e)=\bm g^TA^{-i}\bm g,\quad i=3,4.
$$
Thus, if $|\e|$ sufficiently small, we have
$$
3d\e^2\cdot\frac{\bm g(\e)^T(A(\e)+\vartheta_{\e} I)^{-4}\bm g(\e)}{\left(\bm g(\e)^TA(\e)^{-3}\bm g(\e)\right)^2}\leq 3d\e^2\cdot\frac{\bm g(\e)^TA(\e)^{-4}\bm g(\e)}{\big(\bm g(\e)^TA(\e)^{-3}\bm g(\e)\big)^2}\leq 1,
$$
and it follows from  \eqref{eq2018} and \eqref{eq1416} that
$$
    p_{\e}(\xi\e^2)\leq\|\widetilde{\bm x}\|^2 -d\e^2\leq  \D^2\leq \|\widetilde{\bm x}\|^2= p_{\e}(0).
$$

By monotonicity and continuity of $p_{\e}(\la)$, we know that there is a unique $\hat{\la}\in[0,~\xi\e^2]$, such that $p_{\e}(\hat{\la})=\big\|(A(\e)+\hat{\la} I)^{-1}\bm g(\e)\big\|^2=\D^2$. Let
$\widehat{\bm x}=-(A(\e)+\hat{\la} I)^{-1}\bm g(\e)$, then we have that
$$
\hat{\la}\geq 0,~~(A(\e)+\hat{\la} I)\widehat{\bm x}=-\bm g(\e),~~\hat{\la}(\D-\|\widehat{\bm x}\|)=0 ~~{\rm and}~~A(\e)+\hat{\la} I\succ\bm O.
$$
When $|\e|$ is sufficiently small, it is seen from Theorem \ref{7.59} that $(\hat{\la},\widehat{\bm x})=(\la_*(\e),\bm x_*(\e))$ is the Lagrange multiplier and optimal solution of TRS \eqref{2}, respectively.
So we have
$$
0\leq\lim_{\e\rightarrow 0}\la_*(\e)\leq\lim_{\e\rightarrow 0}\xi\e^2=0,~~{\rm and}~~0\leq
\lim_{\e\rightarrow 0}\frac{\la_*(\e)-\la_*}{\e}=\lim_{\e\rightarrow 0}\frac{\la_*(\e)}{\e}\leq\lim_{\e\rightarrow 0}\xi\e=0,
$$
where we used $\la_*=0$.
That is, $\underset{\e\rightarrow 0}{\lim}\la_*(\e)=0=\la_*$ and $\la'(0)=0$. It then
follows from \eqref{eq16.24} and \eqref{eq1659} that
$$
\lim_{\e\rightarrow 0}\bm x_*(\e)=\bm x_*~~{\rm and}~~
   \lim_{\e\rightarrow 0}\frac{\bm x_*(\e)-\bm x_*}{\e}=   -A^{-1}(E\bm x_*+\bm e),
$$
which completes the proof.
\end{proof}
\begin{remark}
In  summary, if ${\la_{*}>0}$ or ${\la_{*}=0}$ and ${\|\bm x_{*}\|<\D}$, both $\la_*(\e)$ and $\bm x_*(\e)$ are continuous and differentiable at $\e=0$. However, when
${\la_{*}=0}$ and ${\|\bm x_{*}\|=\D}$, if $\bm x^T_*A^{-1}(\bm e+E\bm x_*)\neq 0$, then $\la_*(\e)$ and $\bm x_*(\e)$ are only continuous but are not differentiable at $\e=0$.
\end{remark}

\section{Numerical Experiments}

In this section, we show that our theoretical results are computable for large-scale problems. We discuss how to evaluate the scalars $\eta_1$ and $\eta_2$ introduced in Theorem \ref{thm2119}, as well as the condition numbers $s(\la_*)$ and $s(\bm x_*)$ for the Lagrange multiplier $\la_*$ and the solution $\bm x_*$ of the TRS \eqref{1}. All the numerical experiments were run on a AMD R7 5800H  CPU 3.20 GHz with 16GB RAM under Windows 11 operation system. The experimental results are obtained from using MATLAB R2022a implementation with machine precision ${u}_{\rm mach}\approx2.22\times 10^{-16}$.

As $A$ is a symmetric matrix, we use the $k$-step Lanczos process to generate an orthonormal basis $Q_k=[\bm q_1,\bm q_2,\ldots,\bm q_{k}]$ for the Krylov subspace \cite{2,Y.S}
$$
\mathcal{K}_k(A,\bm g)=span\{\bm g, A\bm g,\ldots,A^{k-1}\bm g\},
$$
where
$\bm q_1=\bm g/\|\bm g\|$. The following Lanczos relation holds
\begin{equation}\label{21.18}
 \begin{aligned}
 AQ_k&=Q_kT_k+\theta_{k}\cdot\bm q_{k+1}\bm e_k^T,
 \end{aligned}
\end{equation}
 where $\bm e_k$ is the $k$-th column of the identity matrix, and
 $$
 T_k =Q^T_kAQ_k=\begin{pmatrix}
 \varpi_1 & \theta_1 & &\\
 \theta_1  &  \varpi_2& \theta_2 & \\
          &  \ddots &  \ddots& \ddots&    \\
          &         & \theta_{k-2}& \varpi_{k-1} & \theta_{k-1}\\
          &         &            &  \theta_{k-1} &  \varpi_{k}
 \end{pmatrix}
 \in \mathbb{R}^{k\times k}
$$
is a tridiagonal matrix.

First, we consider how to evaluate $\eta_1$ and $\eta_2$ defined in Theorem \ref{thm2119}, by using the Lanczos method.
Let $(\widehat{\alpha}_k^{(k)},\bm z_k)$ be the smallest eigenpair and $\widehat{\alpha}_{k-1}^{(k)}$ be the second smallest eigenvalue of $T_k$, respectively.
Without loss of generality, suppose that $\bm g$ is not orthogonal to the eigenspace corresponding to $\alpha_{n-1}$ and $\alpha_n$.
Inspired by the results given in \cite{SS}, as $k$ increases, we make use of $\widehat{\alpha}_k^{(k)},\widehat{\alpha}_{k-1}^{(k)}$ as approximations to $\alpha_n$ and $\alpha_{n-s}$, respectively,
and exploit $Q_k\bm z_k$ to approximate ${U_2U_2^T\bm g}/{\|U_2U_2^T\bm g\|}$.
Consequently, we can use $\bm e_1^T\bm z_k={\bm g^TQ_k\bm z_k}/{\|\bm g\|}$ to approximate
$\cos\angle(\bm g,\mathcal{U}_2)=\frac{\|U_2^T\bm g\|}{\|\bm g\|}$ as $k$ increases.
Therefore,
\begin{align}
\eta^{(k)}_1&=\frac{1}{2}\min\Bigg\{
\frac{(\widehat{\alpha}_{k-1}^{(k)}-\widehat{\alpha}_k^{(k)})
\big(\|\bm g\|\big\|\big(T_k-\widehat{\alpha}_k^{(k)}I\big)^{\dag}\bm e_1\big\|-\D\big)}{\frac{3}{2}\|\bm g\|\big\|\big(T_k-\widehat{\alpha}_k^{(k)}I\big)^{\dag}\bm e_1\big\|+1+\|\bm g\|\cdot\frac{\bm e_1^T\bm z_k}{\widehat{\alpha}_{k-1}^{(k)}-\widehat{\alpha}_k^{(k)} }},~\widehat{\alpha}_{k-1}^{(k)}-\widehat{\alpha}_k^{(k)}\Bigg\},
\end{align}
and
\begin{align}
\eta^{(k)}_2&=\frac{1}{2}\!~\!\min\!
\left\{\frac{\|\bm g\|\big(\widehat{\alpha}_{k-1}^{(k)}\!-\!\widehat{\alpha}_k^{(k)}\big)
\bm e_1^T\bm z_k}{4\|\bm g\|+\widehat{\alpha}_{k-1}^{(k)}\!-\!\widehat{\alpha}_k^{(k)}},~\frac{\|\bm g\|\bm e_1^T\bm z_k}{2\D},
~\widehat{\alpha}_{k-1}^{(k)}\!-\!\widehat{\alpha}_k^{(k)}\right\},
\end{align}
can be utilized as approximations to $\eta_1$ and $\eta_2$, respectively.

Second, we consider how to compute the condition numbers defined in \eqref{eq3.9} and \eqref{eq3.16}.
Recall that
$$
s(\la_*)= \frac{\|A_*^{-1}\bm x_*\|}{\bm x_*^TA_*^{-1}\bm x_*}\max\left\{1,{\D}\right\}~~{\rm and}~~s(\bm x_*)=\|A_*\|\|(PA_*P)^\dag\|
~~{\rm with}~~P=I-\frac{\bm x_*\bm x_*^T}{\D^2},
$$
respectively. The generalized Lanczos trust-region (GLTR) method is a popular approach for solving the large-scale TRS \eqref{1} \cite{10}. Indeed, it solves the following problem instead
 \begin{equation}\label{9.00}
 \begin{aligned}
 \min_{\|\bm x\|\leq\D,~ \bm x\in \mathcal{K}_k(A,\bm g)}
 \left\{f(\bm x)=\frac{1}{2}\bm x^T A\bm x+\bm x^T\bm g\right\},
 \end{aligned}
 \end{equation}
 which reduces to a TRS of size $k$-by-$k$
 \begin{equation}\label{9.25}
 \begin{aligned}
   \min_{ \|\bm h\|\leq \D}\left\{f_k(\bm h)= \frac{1}{2}\bm h^T T_k\bm h+\|\bm g\|\cdot\bm h^T \bm e_1\right\}.
   \end{aligned}
 \end{equation}
If we denote by
 $
 \bm h_k=\arg\min_{\|\bm h\|\leq \D}f_k(\bm h),
 $
then $\bm x_k=Q_k\bm h_k\in \mathcal{K}_k(A,\bm g)$ is the minimizer of  \eqref{9.00}, which can be used as
an approximation to $\bm x_{*}$. Let $\la_k$ be a Lagrangian multiplier of
TRS \eqref{9.25}. Then it follows that \cite{BF,5,4}
\begin{equation*}
\begin{aligned}
\|\bm x_k-\bm x_*\|\leq\mathcal{O}\bigg( \Big(   \frac{\sqrt{\kappa_*}-1}{\sqrt{\kappa_*}+1}\Big)^k  \bigg)~~{\rm  and}~~ |\la_k-\la_*|\leq \mathcal{O}\bigg(k \Big(   \frac{\sqrt{\kappa_*}-1}{\sqrt{\kappa_*}+1}\Big)^{2k}  \bigg),
\end{aligned}
\end{equation*}
where $\kappa_*=\|A_*\|\|A_*^{-1}\|$.
Let $\widehat{T}_k=T_k+\la_k I$, then $\widehat{T}_k\succ\bm O$ \cite[Theorem 5.3]{10}. As a result, we use
$$
s(\la_k)= \frac{\|\widehat{T}_k^{-1}\bm h_k\|}{\bm h_k^T\widehat{T}_k^{-1}\bm h_k}\max\left\{1,{\D}\right\}~~{\rm and}~~s(\bm h_k)=\|\widehat{T}_k\|\|(P_k\widehat{T}_kP_k)^\dag\|,
$$
with $P_k=I_k-\frac{\bm h_k\bm h_k^T}{\D^2}$ to approximate $s(\la_*)$ and $s(\bm x_*)$, respectively.
Notice that $P_k\widehat{T}_kP_k\succeq\bm O$. By Theorem \ref{thm111} (ii), $\|(P_k\widehat{T}_kP_k)^\dag\|$ is equal to the reciprocal of the $(k-1)$-th largest eigenvalue of $P_k\widehat{T}_kP_k$.

\begin{table}[t]
\begin{center}\caption{``Exact values" of the parameters}\label{tab:21.3444}
\begin{tabular}{c|c|c|c|c|c}
\hline
 &\!\!\!$\la_* $ \!\!\!   &  \!\!\!$s(\la_*) $ \!\!\!   &  \!\!\! $s(\bm x_{*})$&$\eta_1$ &  $\eta_2$ \\
\hline
{\small Problem (a)}\!\!&\!\! 2.5004e+03\!\! & \!\!1.1284 \!\! & \!\!\!  2.3448e+03\!\! &3.600e$-$03& 1.7615e$-$03 \\
   \hline
{\small Problem (b)}\!\!&\!\! 7.5011e+03\!\!&\!\!1.0723\!\!& \!\!\!  1.7313e+03\!\! &$-$0.1455& 3.5104e$-$03  \\
   \hline
{\small Problem (c)}\!\!& \!\! 4.3949\!\!
 &   \!\!1.1496\!\!& \!\!\!   182.9145\!\!&0.0120&1.0348e-04\\
   \hline
{\small Problem (d)}\!\!&\!\! 6.7513\!\!
 &  \!\! 1.1268\!\!& \!\!\!   158.9863\!\! &\!0.0145\!  &\! 1.0348e-04 \\
   \hline
\end{tabular}
\end{center}
\end{table}

\begin{figure}[H]\caption{Convergence curves the four parameters during Lanczos iterations}\label{fig:1047}
	\begin{subfigure}{0.5\textwidth}
		\centering
		\includegraphics[width=6cm,height=4cm]{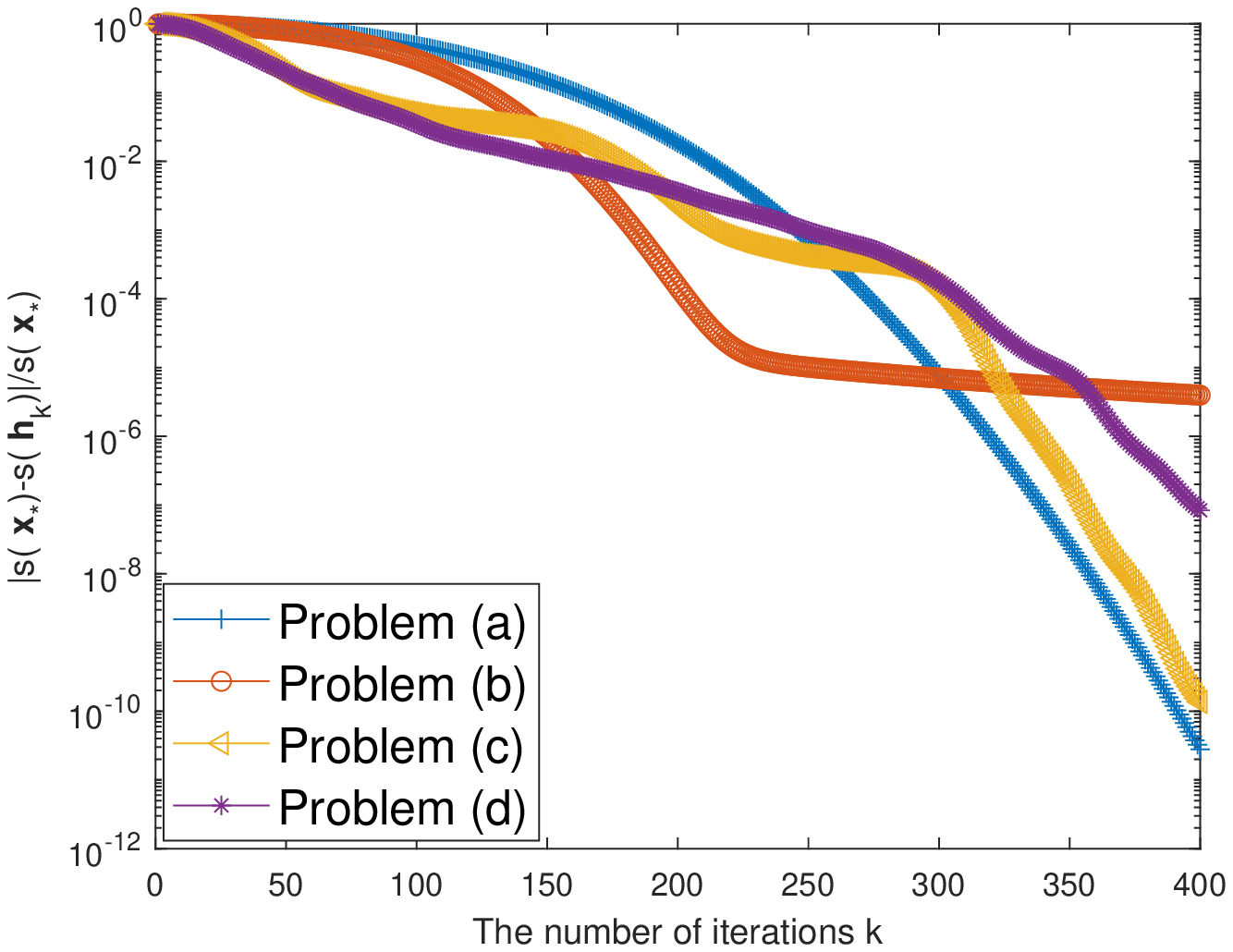}
		\caption*{(i) Convergence  curves of $\frac{|s(\bm h_k)-s(\bm x_*)|}{s(\bm x_*)}$}         
	\end{subfigure}
	\begin{subfigure}{0.5\textwidth}
		\centering
		\includegraphics[width=6cm,height=4cm]{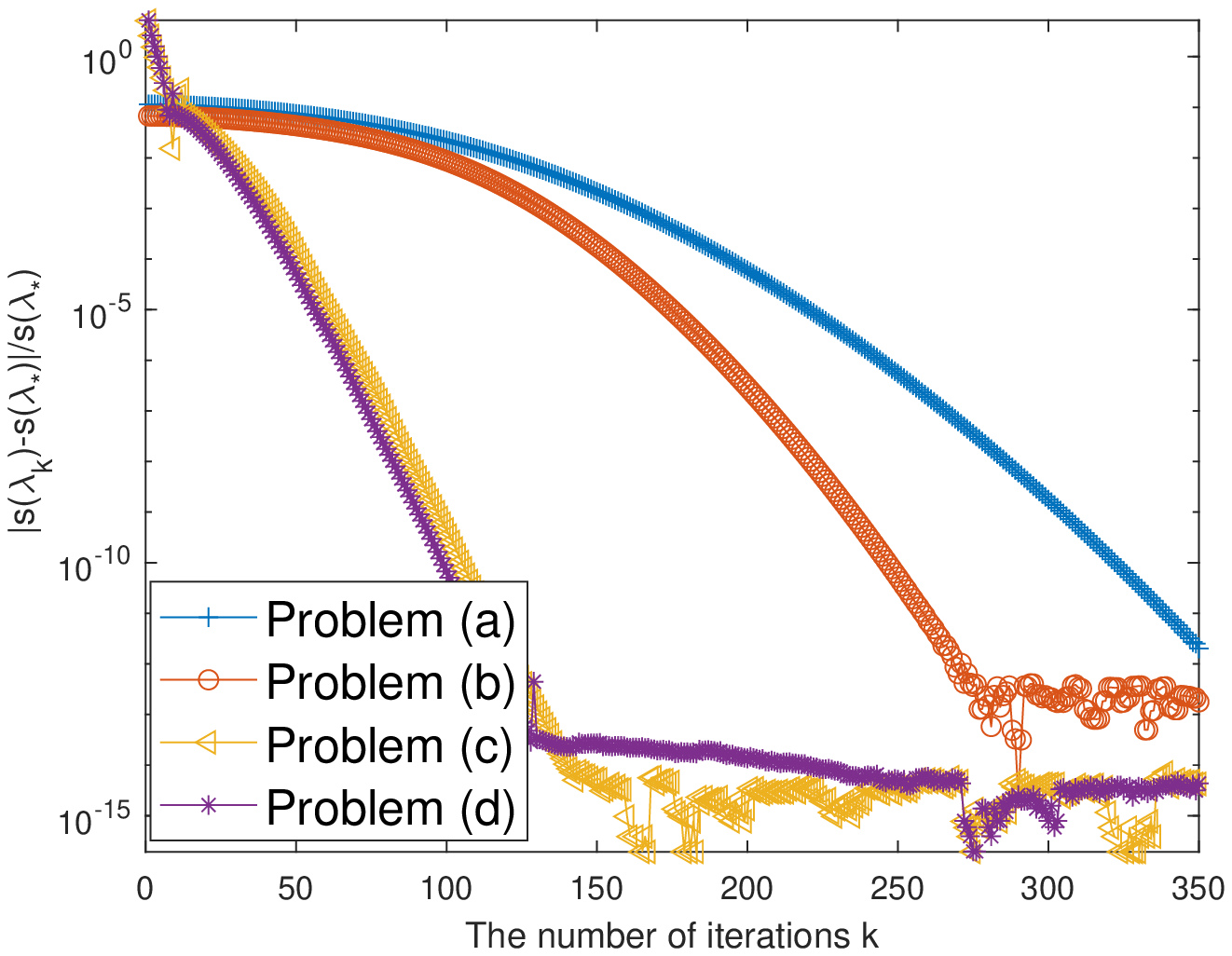}
		\caption*{(ii) Convergence curves of $\frac{|s(\la_k)-s(\la_*)|}{s( \la_*)}$}         
	\end{subfigure}
\begin{subfigure}{0.5\textwidth}
		\centering
		\includegraphics[width=6cm,height=4cm]{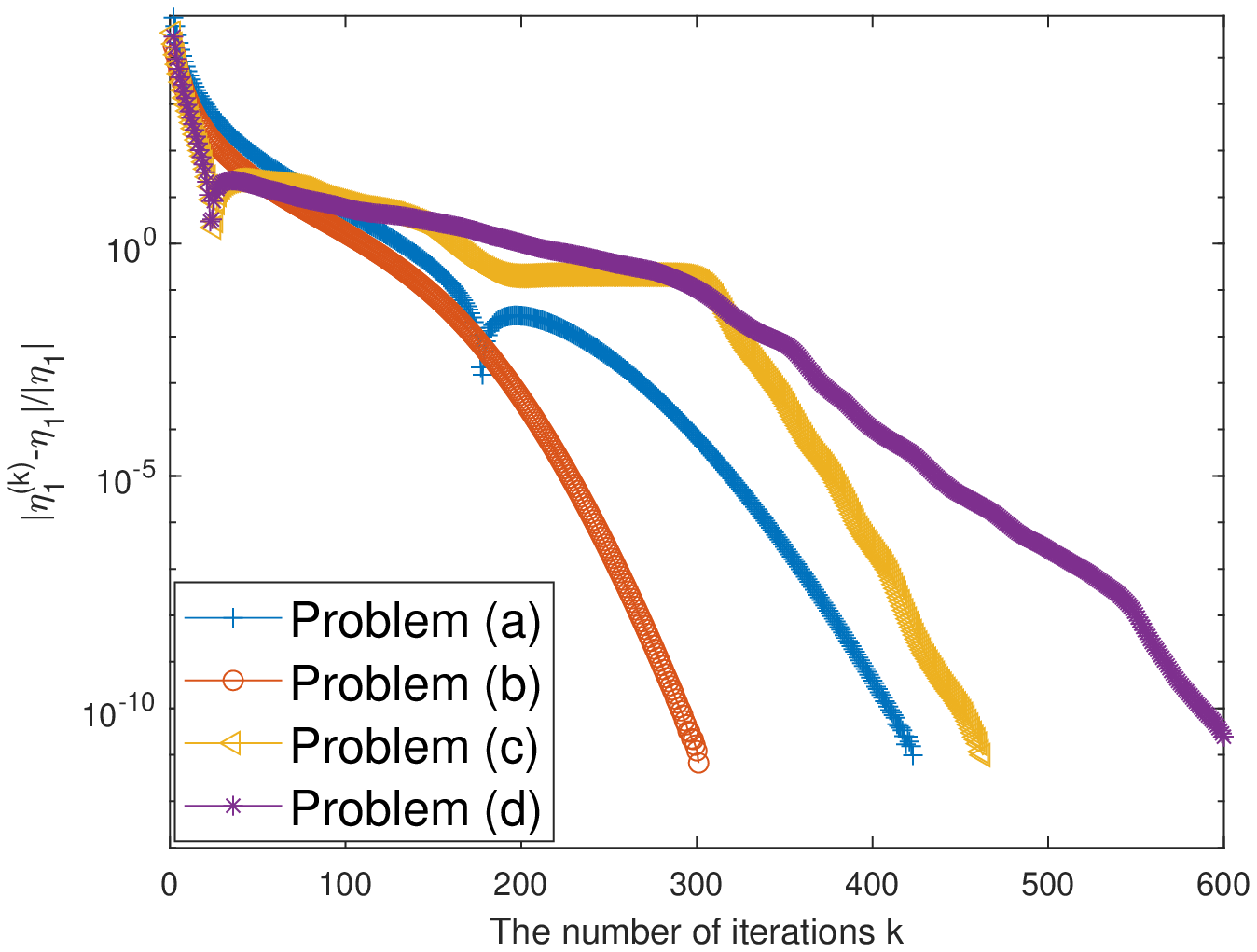}
		\caption*{(iii) Convergence curves of $\frac{|\eta^{(k)}_1-\eta_1|}{|\eta_1|}$}         
	\end{subfigure}
\begin{subfigure}{0.5\textwidth}
		\centering
		\includegraphics[width=6cm,height=4cm]{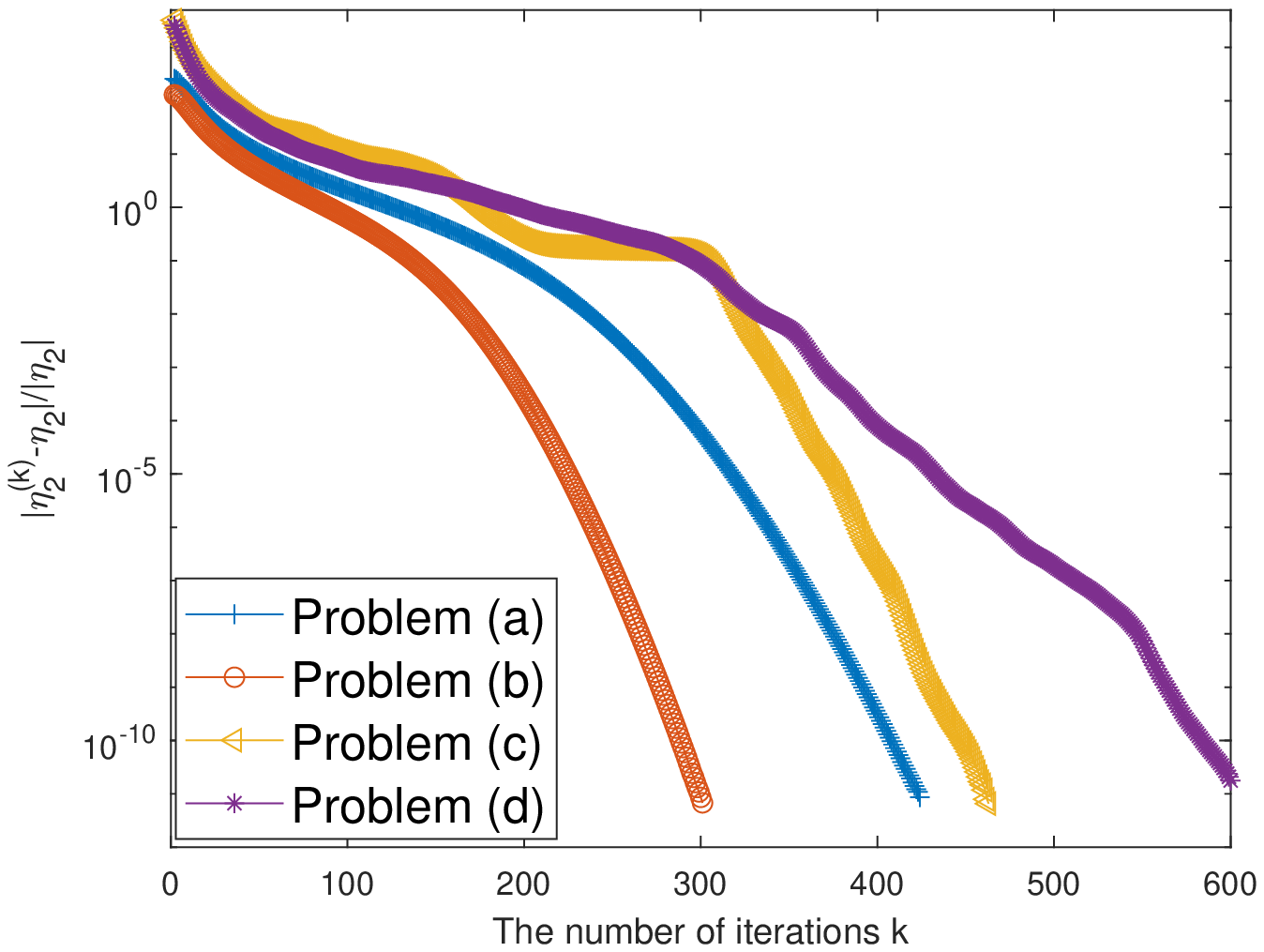}
		\caption*{(iv) Convergence curves of $\frac{|\eta^{(k)}_2-\eta_2|}{|\eta_2|}$}        
	\end{subfigure}
\end{figure}

To illustrate the efficiency of our strategy, we consider the TRS \eqref{1} with
$$
A = diag(\alpha_1,\alpha_2,\ldots,\alpha_n)\in \mathbb{R}^{n\times n}~,~\bm g = (1,1,\ldots,1)^T\in\mathbb{R}^n~~{\rm and}~~\D=1,
$$
where $n=5000$ and the $\{\alpha_i\}$'s are chosen in the following four ways
\begin{align*}
&{\rm Problem~  (a)}:~\alpha_i=i-\frac{n}{2},~~i=1,2,\ldots,n; \\
&{\rm Problem~  (b)}:~\alpha_i=\frac{i^2}{n}-\frac{n}{2}+\frac{1}{n},~~i=1,2,\ldots,n;\\
&{\rm Problem~  (c)}:~(\alpha_1,\alpha_2,\ldots,\alpha_n)={\tt unifrnd(1,1000,1,n)};\\
&{\rm Problem~  (d)}:~(\alpha_1,\alpha_2,\ldots,\alpha_n)={\tt unifrnd(-1,1000,1,n)}.\\
\end{align*}
Here Problems (a) and (b) were used in \cite{G.S.}, and $\tt unifrnd(a,b,m,n)$ is the MATLAB build-in function for generating an $m\times n$ random matrix, whose elements are in  continuous uniform distribution, with the lower endpoint $ a$ and the upper endpoint $ b$.

In this example, we make use of the MATLAB build-in
function $\tt trust.m$ to compute
$\la_*$ and $\bm x_*$, and
$$
\frac{\|(A+\la_*I)\bm x_*+\bm g\|}{\|\bm g\|}\leq 5\times 10^{-16}~~ {\rm and} ~~\frac{|\|\bm x_*\|-\D|}{\D}\leq 5\times 10^{-14}.
$$
The numerical results are summarized in Table \ref{tab:21.3444} and Figure \ref{fig:1047}, respectively. It is obvious to see that these parameters can be approximated very well.
On the other hand, it is seen from Figure \ref{fig:1047} (i) that the convergence of $s(\bm x_*)$ may be (relatively) slow. Indeed, it is only necessary to estimate the order of the condition number, and a rough evaluation is enough in practice.

\section{Conclusion}
Trust-region subproblem plays an important role in the areas of numerical linear algebra and numerical optimization. In this paper, we focus on the perturbation analysis of TRS and establish some first-order perturbation results on this problem. First, if TRS \eqref{1} is in \emph{easy case}, we derive an upper bound on $|\e|$, such that the perturbed TRS \eqref{2} is also in \emph{easy case}. Second, we prove that $\la_*(\e)$ and $\bm x_*(\e)$ are continuous at $\e=0$, and $\la_*(0)=\la_*$, $\bm x_*(0)=\bm x_*$. Third, we establish some first-order perturbation bounds and define condition numbers on $\la_*$ and $\bm x_*$, respectively. Specifically,  we point out that as $\la_*>0$, the condition number of the $\bm x_*$ is not necessarily large even if the TRS \eqref{1} is in \emph{nearly hard case}.
Fourth, we indicate that when ${\la_{*}=0}$ and ${\|\bm x_{*}\|=\D}$, if $\bm x^T_*A^{-1}(\bm e+E\bm x_*)\neq 0$, then $\la_*(\e)$ and $\bm x_*(\e)$ are only continuous but are not differentiable at $\e=0$.

The established results are useful to evaluate the ill-conditioning of a TRS problem beforehand, and to access the quality of the approximate solution.
Furthermore, we present some strategies to compute these results by using the Lanczos method for large-scale TRS. Numerical experiments demonstrate that the proposed schemes are effective. On the other hand, there are some problems need to investigate further. For instance, we only consider the situations that TRS is in {\it easy case} or in {\it nearly hard case}. Perturbation analysis on {\it hard case} is much more complicated and deserves further study.

\end{document}